\documentclass[11pt,english]{article}

\usepackage[margin = 2.5cm]{geometry}

\usepackage{amsthm}
\usepackage{amsmath}
\usepackage{amssymb}
\usepackage{mathtools}
\usepackage{graphicx}
\usepackage[hidelinks]{hyperref}
\usepackage{enumitem}

\usepackage{caption}
\usepackage[square,sort,comma,numbers]{natbib} 
\usepackage{setspace}
\usepackage{cleveref}
\theoremstyle{plain}

\newtheorem*{thm*}{Theorem}
\newtheorem{thm}{Theorem}[section]
\Crefname{thm}{Theorem}{Theorems}

\newtheorem*{lem*}{Lemma}
\newtheorem{lem}[thm]{Lemma}
\Crefname{lem}{Lemma}{Lemmas}

\newtheorem*{claim*}{Claim}
\newtheorem{claim}[thm]{Claim}
\crefname{claim}{Claim}{Claims}
\Crefname{claim}{Claim}{Claims}

\newtheorem{prop}[thm]{Proposition}
\Crefname{prop}{Proposition}{Propositions}

\newtheorem{cor}[thm]{Corollary}
\crefname{cor}{Corollary}{Corollaries}

\newtheorem{conj}[thm]{Conjecture}
\crefname{conj}{Conjecture}{Conjectures}

\newtheorem{qn}[thm]{Question}
\Crefname{qn}{Question}{Questions}

\newtheorem{obs}[thm]{Observation}
\Crefname{obs}{Observation}{Observations}

\theoremstyle{definition}

\Crefname{prob}{Problem}{Problems}

\Crefname{defn}{Definition}{Definitions}

\theoremstyle{remark}

\theoremstyle{definition}

\newtheorem{ex}[thm]{Example}
\Crefname{ex}{Example}{Examples}

\usepackage{xpatch}
\xpatchcmd{\proof}{\itshape}{\normalfont\proofnamefont}{}{}
\newcommand{\proofnamefont}{}
\renewcommand{\proofnamefont}{\bfseries}

\usepackage{framed}

\usepackage{soul}
\usepackage{color}

\newcommand{\remove}[1]{}

\newcommand{\ceil}[1]{
    \left\lceil #1 \right\rceil
}
\newcommand{\floor}[1]{
    \left\lfloor #1 \right\rfloor
}

\newcommand{\T}{\mathcal{T}}

\newcommand{\B}{\mathcal{B}}
\newcommand{\C}{\mathcal{C}}
\renewcommand{\Pr}{\mathbb{P}}
\newcommand{\LL}{\mathcal{L}}

\newcommand{\eps}{\varepsilon}

\DeclareMathOperator{\bip}{bip}

\DeclareMathOperator{\pack}{pack}
\newcommand{\fpack}{\pack}

\usepackage{bbm}
\newcommand{\one}{\mathbbm{1}}

\usepackage{multirow}
\usepackage[table]{xcolor}
\def \cola{yellow}
\def \colb{green!90!black}
\def \colc{red!70}
\def \cold{blue!50}

\title{Monochromatic triangle packings in red-blue graphs}

\author{
	Vytautas Gruslys\thanks{
		Email: \texttt{vytautas.gruslys}@\texttt{gmail.com}.
	}
	\and
    Shoham Letzter\thanks{
		Department of Mathematics, 
		University College London, 
		Gower Street, London WC1E~6BT, UK. 
		Email: \texttt{s.letzter}@\texttt{ucl.ac.uk}. 
		Research supported by the Royal Society.
    }
}

\begin{document}

\date{}
\maketitle

\begin{abstract}

	We prove that in every $2$-edge-colouring of $K_n$ there is a collection of $n^2/12 + o(n^2)$ edge-disjoint monochromatic triangles, thus confirming a conjecture of Erd\H{o}s. We also prove a corresponding stability result, showing that $2$-colourings that are close to attaining the aforementioned bound have a colour class which is close to bipartite.
	As part of our proof, we confirm a recent conjecture of Tyomkyn about the fractional version of this problem.

	\setlength{\parskip}{\medskipamount}
    \setlength{\parindent}{0pt}
    \noindent

\end{abstract}
	\section{Introduction} \label{sec:intro}

		A result of Goodman \cite{goodman} shows that every $2$-edge-colouring of $K_n$ contains at least $n^3/24 + o(n^3)$ monochromatic triangles, an estimate which can be seen to be asymptotically tight by letting one of the colour classes be a balanced complete bipartite graph $K_{\ceil{n/2}, \floor{n/2}}$. Erd\H{o}s \cite{erdos} considered\footnote{Erd\H{o}s \cite{erdos} attributes the question to Ordman, Faudree and himself, but in subsequent publications, including \cite{erdos-et-al} which is coauthored by Erd\H{o}s and Faudree, the problem is attributed only to Erd\H{o}s.} a variant of Goodman's result: how many \emph{pairwise edge-disjoint} monochromatic triangles are there guaranteed to be in a $2$-colouring of $K_n$? Again considering the example where one of the colour classes is the balanced complete bipartite graph, Erd\H{o}s made the following conjecture.

		\begin{conj}[Problem 14 in \cite{erdos}] \label{conj:erdos}
			Every $2$-coloured $K_n$ contains $n^2/12 + o(n^2)$ pairwise edge-disjoint monochromatic triangles.
		\end{conj}

		The first progress towards this conjecture was made by Erd\H{o}s, Faudree, Gould, Jacobson and Lehel \cite{erdos-et-al}, who proved that there are always at least $3n^2/55 + o(n^2)$ edge-disjoint monochromatic triangles. To prove this bound they calculated the minimum possible number of edge-disjoint monochromatic triangles in a $2$-coloured $K_{11}$, and used Wilson's theorem which guarantees the existence of an almost decomposition of the edges of $K_n$ into copies of $K_{11}$. This was improved by Keevash and Sudakov \cite{keevash-sudakov}, who proved a bound of $n^2/12.89 + o(n^2)$. They used a reduction to a fractional version of the problem due to Haxell and R\"odl \cite{haxell-rodl}, averaging arguments to relate the answers for $n$ and $n-1$ and also for $3n$ and $n$, and a computer search to calculate the optimal value for $n = 15$.

		Alon and Linial (see \cite{keevash-sudakov}) suggested to study a weaker version of \Cref{conj:erdos}, where the only $2$-colourings allowed are those with a triangle-free colour class. In other words, is it true that in every $n$-vertex graph, whose complement is triangle-free, there are $n^2/12 + o(n^2)$ edge-disjoint triangles? Yuster \cite{yuster07} considered this question and showed that any counterexamples (if exist) have between $0.2501n^2$ and $3n^2/8$ edges. Recently, Tyomkyn \cite{tyomkyn} answered Alon and Linial's question affirmatively. Like Keevash and Sudakov \cite{keevash-sudakov}, he used the reduction of Haxell and R\"odl \cite{haxell-rodl} to the fractional version of the problem. His proof is inductive, using an averaging argument also used by Keevash and Sudakov, and a computer search to resolve the question for small values of $n$. A key ingredient in his argument is a lemma that asserts that graphs which are `critical', namely their complement is not bipartite but can be made bipartite by the removal of one vertex, have large `fractional triangle packings'.

		Our main result in this paper confirms \Cref{conj:erdos}.

		\begin{thm} \label{thm:main}
			Every $2$-coloured $K_n$ contains a collection of $n^2/12 + o(n^2)$ pairwise edge-disjoint monochromatic triangles.
		\end{thm}

		As in \cite{keevash-sudakov,tyomkyn}, it suffices to prove a fractional analogue of \Cref{thm:main}. Our proof of the fractional version is inductive, with a computer search to deal with small value of $n$. Our main inductive step deals with $2$-colourings of $K_n$ where one of the colours is close to bipartite. In particular, it allows us to prove a conjecture of Tyomkyn \cite{tyomkyn} about $2$-coloured complete graphs whose `monochromatic fractional triangle packing number' is close to extremal (see \Cref{thm:frac-main-almost-extremal}; we introduce the relevant notions in the next section). 
		However, for smaller value of $n$, which are too large for the computer search, we also need to consider colourings that are close to `pentagon blow-ups'. An important ingredient in our proof of the almost bipartite case is a result about fractional triangle packings in almost complete graphs (see \Cref{thm:n-four}, which we prove in a separate paper \cite{us2}, in order to keep this paper from becoming unreasonably long). 
		 
		We say that a graph $G$ is \emph{$k$-close to bipartite} if $G$ can be made bipartite by removing at most $k$ edges. If $G$ is not $k$-close to bipartite, we say that it is \emph{$k$-far from bipartite}.
		Tyomkyn \cite{tyomkyn} proved that $n$-vertex graphs, whose complement is triangle-free, and which do not have significantly more than $n^2/4$ edge-disjoint triangles, are close to bipartite. We generalise his result, thus obtaining a stability version of \Cref{thm:main}.

		\begin{thm} \label{thm:main-stability}
			For every $\eps > 0$ there exists $\delta > 0$ such that the following holds for every sufficiently large $n$. If $G$ is a $2$-colouring of $K_n$ where both colour classes are $\eps n^2$-far from bipartite, then there is a collection of $n^2/12 + \delta n^2$ edge-disjoint monochromatic triangles in $G$.
		\end{thm}

		Following Tyomkyn \cite{tyomkyn}, we use a result of Alon, Shapira and Sudakov \cite{alon-shapira-sudakov} (see \Cref{thm:a-s-s}) about the structure of graphs that are far from a given monotone graph property. Additionally, we use a stability version of a fractional version of our main result. 

		In the next section we give an overview of our proof and of the structure of the paper.

	\section{Overview} \label{sec:overview}

		\subsection{A reduction to fractional triangle packings}

			A \emph{triangle packing} in a graph $G$ is a collection of edge-disjoint triangles. 
			Given a graph $G$, we write $\nu(G)$ to denote the number of triangles in a largest triangle packing in $G$. 
			Let $\T(G)$ denote the collection of triangles in a graph $G$.
			A \emph{fractional triangle packing} in a graph $G$ is a function $\omega: \T(G) \to [0,1]$ such that $\sum_{T \in \T(G): \,\, e \subseteq T} \omega(T) \le 1$ for every edge $e$. Let $\nu^*(G)$ denote the weight of the `largest' fractional triangle packing in $G$, namely 
			\begin{equation*}
				\nu^*(G) = \max\left\{\sum_{T \in \T(G)}\omega(T) : \, \text{ $\omega$ is a fractional triangle packing in $G$}\right\}.
			\end{equation*}
			We use the following result of Haxell and R\"odl \cite{haxell-rodl}, which implies that the weight of the largest fractional triangle packing in a graph $G$ is a good approximation for the number of triangles in a largest triangle packing in $G$. 

			\begin{thm}[A special case of Thereom 1 in \cite{haxell-rodl}] \label{thm:haxell-rodl}
				Let $G$ be a graph on $n$ vertices. Then $\nu^*(G) = \nu(G) + o(n^2)$.
			\end{thm}

			Given a fractional triangle packing $\omega$ in a graph $G$ and an edge $e$, we define the weight of $e$, denoted $\omega(e)$, to be the sum of weights of triangles containing $e$, i.e.\ $\omega(e) = \sum_{T \in \T(G): \,\, e \subseteq T} \omega(T)$ (so $\omega(e) \in [0,1]$). 
			We define the \emph{size} of $\omega$, denoted $\omega(G)$, to be the sum of weights of the edges $G$, namely $\omega(G) = \sum_{e \in E(G)} \omega(e)$. Equivalently, $\omega(G) = 3\sum_{T \in \T(G)} \omega(T)$. We use this somewhat non-standard scaling, as it is often more instructive to consider the weights of the \emph{edges} covered by a fractional triangle packing, rather than the weights of the triangles.
			
			Given a red-blue coloured graph $G$, denote by $G_R$ and $G_B$ the subgraphs of $G$ spanned by the red and blue edges, respectively. We define $\fpack(G)$ to be the size of the largest monochromatic triangle packing in $G$, namely
			\begin{align} \label{eq:defn-pack}
				\begin{split}
					\fpack(G) = \, 
					& \max\{\omega(G_R) : \text{$\omega$ is a fractional triangle packing in $G_R$}\} \,\,+ \\
					& \max\{\omega(G_B) : \text{$\omega$ is a fractional triangle packing in $G_B$}\}.
				\end{split}
			\end{align}
			Equivalently, $\fpack(G) = 3(\nu^*(G_R) + \nu^*(G_B))$. The following corollary is an immediate consequence of \Cref{thm:haxell-rodl}.

			\begin{cor} \label{cor:haxell-rodl}
				Let $G$ be a red-blue colouring of $K_n$. Then there is a monochromatic triangle packing in $G$ that covers at least $\fpack(G) + o(n^2)$ edges.
			\end{cor}

		\subsection{Minimising the size of a monochromatic fractional triangle packing}

			By \Cref{cor:haxell-rodl}, in order to tackle \Cref{thm:main}, it suffices to solve the following extremal question: what is the minimum of $\fpack(G)$ among all red-blue colourings of $K_n$? We answer this question as follows.

			\begin{thm} \label{thm:frac-main-extremal}
				Let $n \ge 26$. Suppose that $G$ is a red-blue colouring of $K_n$. Then $\fpack(G) \ge \floor{(n-1)^2/4}$, with equality if and only if one of $G_R$ and $G_B$ is the union of a balanced complete bipartite graph $K_{\ceil{n/2}, \floor{n/2}}$ with a matching.
			\end{thm}

			We note that the statement of \Cref{thm:frac-main-extremal} does not hold for all values of $n$. To see this consider the following example.

			\begin{ex} \label{ex:pentagon-blowup}
				A \emph{pentagon blow-up} is a red-blue colouring of $K_n$ where the vertices can be partitioned into five non-empty sets $A_1, \ldots, A_5$ such that the edges between $A_i$ and $A_{i+1}$ are red, the edges between $A_i$ and $A_{i+2}$ are blue, and the edges in $A_i$ are coloured arbitrarily, for $i \in [5]$ (addition of indices is taken modulo $5$).

				A \emph{balanced pentagon blow-up} is a pentagon blow-up whose blob sizes differ by at most $1$. It is not hard to check that a balanced pentagon blow-up $G$ with blob sizes $a_1, \ldots, a_5$ satisfies $\fpack(G) = 3\sum_i \binom{a_i}{2}$ (see \Cref{prop:pentagon-packing}).
			\end{ex}
			
			In particular, a balanced pentagon blow-up $G$ on $20$ vertices satisfies $\fpack(G) = 3 \cdot 5 \cdot 6 = 90 = \floor{19^2/4}$, so it achieves equality in \Cref{thm:frac-main-extremal}, but neither $G_R$ nor $G_B$ are close to bipartite. Similarly, a balanced pentagon blow-up on $17$ vertices $H$ satisfies $\fpack(H) = 63 < 64 = 16^2/4$, so it violates the inequality in \Cref{thm:frac-main-extremal}.\footnote{Our findings show that the statement of \Cref{thm:frac-main-extremal} holds also for $n \ge 21$, and moreover the inequality holds for $n \ge 18$ (both of these statements are best possible as can be seen by the examples mentioned above), but we do not formally prove this slight generalisation of \Cref{thm:frac-main-extremal} for technical and presentational reasons.}

			Note that our first main theorem, \Cref{thm:main}, follows directly from \Cref{thm:frac-main-extremal} and \Cref{cor:haxell-rodl}. 

		\subsection{Characterising almost extremal examples} \label{subsec:almost-extremal}

			The following observation, due to Keevash and Sudakov \cite{keevash-sudakov}, is a a very useful tool in our arguments; it was also used in \cite{tyomkyn}.

			\begin{obs}[A variant of Lemma 2.1 in \cite{keevash-sudakov}] \label{obs:averaging} 
				Let $G$ be a red-blue colouring of $K_{n+1}$. Then 
				\begin{equation*} 
					\fpack(G) \ge \frac{1}{n-1} \cdot \sum_{u \in V(G)} \fpack(G \setminus \{u\}).
				\end{equation*}
			\end{obs}

			\begin{proof} 
				Let $H$ be a graph on $n+1$ vertices.
				Given a vertex $u$, let $\omega_u$ be a fractional triangle packing of $H \setminus \{u\}$ of maximum size. Let $\omega$ be the fractional triangle packing in $H$ defined by $\omega = \frac{1}{n-1} \sum_{u \in V(H)} \omega_u$. Note that $\omega$ is indeed a fractional triangle packing: given an edge $e = xy$, it receives zero weight from $\omega_x$ and $\omega_y$, and weight at most $1$ from $\omega_u$ for $u \neq x, y$, amounting to weight at most $1$ in $\omega$; in particular, $\omega$ assigns weight at most $1$ to each triangle in $H$. By choice of $\omega_u$, it follows that $\nu^*(H) \ge \frac{1}{n-1} \sum_u \nu^*(H \setminus \{u\})$.
				The required inequality follows by plugging in $H = G_R$ and $H = G_B$, and using that $\fpack(F) = 3(\nu^*(F_R) + \nu^*(F_B))$ for every red-blue coloured graph $F$.
			\end{proof}
			 
			\Cref{obs:averaging} suggests an inductive approach towards determining the minimum of $\fpack(\cdot)$ over red-blue colourings of $K_n$. More precisely, it seems convenient to use the scaling $\frac{\fpack(G)}{n(n-1)}$. Indeed, if $G$ is a red-blue colouring of $K_{n+1}$ with $\fpack(G) \le \alpha \cdot {n(n+1)}$, then \Cref{obs:averaging} implies that for some vertex $u$ we have $\fpack(G \setminus \{u\}) \le \alpha \cdot n(n-1)$. Recall that our goal in \Cref{thm:frac-main-extremal} is to show that the minimum value of $\fpack(\cdot)$ over red-blue colourings of $K_n$ is $\floor{(n-1)^2/4}$. Thus, in order to make use of \Cref{obs:averaging}, one is led to ask: which red-blue colourings $G$ of $K_n$ satisfy $\fpack(G) \le n(n-1)/4$? We answer this question as follows, confirming a conjecture of Tyomkyn \cite{tyomkyn}.

			\begin{thm} \label{thm:frac-main-almost-extremal}
				Let $n \ge 26$. Suppose that $G$ is a red-blue colouring of $K_n$ with $\fpack(G) \le n(n-1)/4$. Then one of $G_R$ and $G_B$ is $(n/8)$-close to bipartite.
			\end{thm}

			We note that the statement in \Cref{thm:frac-main-almost-extremal} does not hold for $n \le 25$. For example, a balanced blow-up of a pentagon on $25$ vertices $H$ has $\fpack(H) = 150 = 25 \cdot 24 / 4$, despite both $H_R$ and $H_B$ being far from bipartite (recall \Cref{ex:pentagon-blowup}).

			\Cref{thm:frac-main-almost-extremal} does not immediately imply \Cref{thm:frac-main-extremal}. Nevertheless, the tools that we develop in order to prove the former theorem allow us to deduce the latter quite easily; see \Cref{sec:frac-extremal}.

			Turning back to the proof of \Cref{thm:frac-main-almost-extremal}, our plan is to prove it by induction. The following lemma provides us with the induction step; we prove it in \Cref{sec:almost-bip}.

			\begin{lem} \label{lem:almost-bip-step}
				Let $n \ge 22$. Let $G$ be a red-blue colouring of $K_{n+1}$, such that $\fpack(G) \le n(n+1)/4$. Suppose that for some vertex $u$ the colouring $H = G \setminus \{u\}$ satisfies $\fpack(H) \le n(n-1)/4$ and $H_B$ is $(n/8)$-close to bipartite. Then $G_B$ is $(n+1)/8$-close to bipartite.
			\end{lem}

			Our approach for dealing with small value of $n$ is via a computer search, again capitalising on \Cref{obs:averaging}. Indeed, by this observation, in order to find all red-blue colourings $G$ of $K_{n+1}$ with $\fpack(G) \le n(n+1)/4$, it suffices to find all red-blue colourings of $H$ of $K_n$ with $\fpack(H) \le n(n-1)/4$, and consider all possible ways of extending $H$ to a red-blue colouring of $K_{n+1}$. 

			Ideally, such a search would tell us that for some $n \ge 22$, if $H$ is a red-blue colouring of $K_{n}$ then one of $H_R$ and $H_B$ is $(n/8)$-close to bipartite. (In fact, due to \Cref{ex:pentagon-blowup}, we would have to take $n \ge 26$.) However, the number of examples makes this unfeasible. Instead, we use the above approach (though with a more sophisticated implementation; see \Cref{sec:computer-search}) to find all relevant examples up to $n = 17$. We then partition the collection of examples into two sets: examples that are close to pentagon blow-ups; and the remaining examples. 
			
			For the former ones -- namely those that are close to pentagon blow-ups -- we show that the only way to extend any of them to a red-blue colouring of a complete graph with small enough value of $\fpack(\cdot)$ is by extending them to a colouring which is again very close to a pentagon blow-up. Iterating this, we show that the `almost pentagon-blow-ups' of order $17$ cannot be extended to examples with small enough $\fpack(\cdot)$ on more than $25$ vertices.

			We keep extending the latter examples -- namely those that are not close to a pentagon blow-up -- until they have $22$ vertices, in which case we find that all surviving examples are close to bipartite, allowing us to apply \Cref{lem:almost-bip-step}.

			The findings of our computer search are summarised in the following lemma; we discuss our algorithm in more detail in \Cref{sec:computer-search}.
			The relevant certificates corresponding to the computer search can be found \href{https://liveuclac-my.sharepoint.com/:f:/g/personal/ucahsle_ucl_ac_uk/EltPNwoeZx5GoXvNBgA3RDcBqGrsLjwjuZ1SW56SrlGK8A?e=5\%3ahICC1Q&at=9}{\textcolor{blue}{here}}.
			\begin{lem} \label{lem:computer}
				Suppose that $G_{0} \subseteq \ldots \subseteq G_{22}$ is a sequence where $G_n$ is a red-blue colouring of $K_n$ with $\fpack(G_n) \le n(n-1)/4$ for $n \in \{0, \ldots, 22\}$. Then one of the following conditions is satisfied.
				\begin{itemize} 
					\item 
						$G_{17}$ is a blow-up of a pentagon with blobs sizes $x_1, \ldots, x_5$ (in some order), where $(x_1, \ldots, x_5) \in \{(3,3,3,4,4), (2,3,4,4,4), (3,3,3,3,5)\}$,
					\item
						$G_{17}$ is one edge-flip away from a pentagon blow-up with blob sizes $3,3,3,4,4$,  
					\item
						the blue edges of $G_{22}$ span a graph which is $2$-close to bipartite,
					\item
						the red edges of $G_{22}$ span a graph which is $2$-close to bipartite.
				\end{itemize}
			\end{lem}
			In the following lemma we take care of the case where one of the first two items in \Cref{lem:computer} hold, namely when $G_{17}$ is close to a pentagon blow-up; we prove it is \Cref{sec:pentagon}.

			\begin{lem} \label{lem:pentagon}
				There is no sequence $G_{17} \subseteq \ldots \subseteq G_{26}$ such that $G_n$ is a red-blue colouring of $K_n$ satisfying $\fpack(G_n) \le n(n-1)/4$ for $n \in \{17, \ldots, 26\}$, and $G_{17}$ is either a pentagon blow-up with blobs sizes $x_1, \ldots, x_5$, where $(x_1, \ldots, x_5) \in  \{(3,3,3,4,4), (2,3,4,4,4), (3,3,3,3,5)\}$; or it is one edge-flip away from a pentagon blow-up with blob sizes $3,3,3,4,4$.
			\end{lem}

			It is now easy to prove \Cref{thm:frac-main-almost-extremal}.

			\begin{proof}[Proof of \Cref{thm:frac-main-almost-extremal}]
				Recall that $G$ is a red-blue colouring of $K_n$ with $\fpack(G) \le n(n-1)/4$, where $n \ge 26$. By \Cref{obs:averaging}, there is a sequence $G_{17} \subseteq \ldots \subseteq G_{n} = G$ such that $G_i$ is a red-blue colouring of $K_i$ with $\fpack(G_i) \le i(i-1)/4$, for $i \in \{17, \ldots, 26\}$.
				By \Cref{lem:computer}, either $G_{17}$ satisfies the conditions of \Cref{lem:pentagon}, or, without loss of generality, the blue edges in $G_{22}$ span a graph which is $2$-close to bipartite. In the former case we reach a contradiction to $n \ge 26$ by \Cref{lem:pentagon}, and in the latter case we conclude that the blue edges in $G_i$ form a graph which is $(i/8)$-close to bipartite, for all $i \in \{22, \ldots, n\}$, as required.
			\end{proof}

		\subsection{Stability}

			As in \cite{tyomkyn}, a result of Alon, Shapira and Sudakov \cite{alon-shapira-sudakov} implies that in order to prove our stability result, \Cref{thm:main-stability}, it suffices to show that in every red-blue colouring $G$ of $K_n$ with $\fpack(G) \le n(n-1)/4 + \eta n$ one of the colour classes is close to bipartite, for some constant $\eta > 0$ and sufficiently large $n$ (see \Cref{sec:stability}). The following theorem provides us with such a result.

			\begin{thm} \label{thm:frac-main-stability}
				There exists $\eta > 0$ such that the following holds for all sufficiently large $n$.
				Let $G$ be a red-blue colouring of $K_n$. Then either $\fpack(G) \ge n(n-1)/4 + 2\eta n$, or one of $G_B$ and $G_R$ is $(1/8 + \eta) n$-close to bipartite.
			\end{thm}

			In a sense, this result seems harder to prove than \Cref{thm:frac-main-almost-extremal}, because here we need to consider colourings with slightly larger value of $\fpack(\cdot)$. However, with \Cref{thm:frac-main-almost-extremal} in hand, we get the induction base for \Cref{thm:frac-main-stability} for free (by taking sufficiently small $\eta$), and it also allows us to start at a larger value of $n$. Unfortunately, we do need to repeat some of the arguments used in the proof of \Cref{thm:frac-main-almost-extremal}, but the larger value of $n$ allows for a simpler presentation.

		\subsection{A key lemma}

			A \emph{triangle decomposition} in a graph $G$ is a collection of edge-disjoint triangles that covers all the edges in $G$. Similarly, a \emph{fractional triangle decomposition} in a graph $G$ is a fractional triangle packing where every edge has weight $1$.

			The following result, which shows that almost complete graphs have fractional triangle decompositions, is key in our arguments regarding potential examples that are close to bipartite. 

			\begin{thm} \label{thm:n-four}
				Let $G$ be a graph on $n \ge 7$ vertices with $e(G) \ge \binom{n}{2} - (n - 4)$. Then there is a fractional triangle decomposition in $G$.
			\end{thm}

			This result is tight in two ways: the complete graph on six vertices with two edges removed (intersecting or not) does not have a fractional triangle decomposition; and the graph on vertex set $[n]$ with non-edges $\{xn : x \in \{4, \ldots, n-1\}\} \cup \{12\}$ is an $n$-vertex graph with $n - 3$ non-edges that does not have a fractional triangle decomposition. 

			A well-known conjecture of Nash-Williams \cite{nash-williams} asserts that every $n$-vertex graph $G$ with minimum degree at least $3n/4$, where $n$ is large and $G$ satisfies certain `divisibility conditions', has a triangle decomposition. While this conjecture is still very much open, significant progress towards it has been made. Recently, Delcourt and Postle \cite{delcourt-postle} showed that every $n$-vertex graph with minimum degree at least $0.83 n$ has a fractional triangle decomposition, improving on several previous results (see, e.g., \cite{gustavsson,yuster05,dukes,garaschuk,dross}).
			Combined with a result of Barber, K\"uhn, Lo and Osthus \cite{barber-et-al}, it follows that the statement obtained by replacing $3/4$ by $0.831$ in Nash-Williams's conjecture holds. 
			This result of Delcourt and Postle (or any result about fractional triangle decompositions in graphs with large minimum degree) can be used to prove \Cref{thm:n-four} for sufficiently large $n$. 
			
			However, crucially, we need \Cref{thm:n-four} to hold for all $n \ge 7$. We thus prove \Cref{thm:n-four} ourselves. Due to the length of the proof and of the current paper, we prove the theorem in a separate paper \cite{us2}. Our proof is again inductive, using an averaging argument as in \Cref{obs:averaging}, with a computer search to prove the base case, but the details are rather involved. 

			The following fractional version of \Cref{thm:n-four} can be deduced from \Cref{thm:n-four} (see \cite{us2}). 

			\begin{cor} \label{cor:n-four-frac}
				Let $G$ be a complete graph on $n \ge 7$ vertices, and let $\phi : E(G) \to [0,1]$ be such that $\sum_{e \in E(G)} \phi(e) \ge \binom{n}{2} - (n-4)$. 
				Then there is a fractional triangle packing $\omega$ in $G$ such that $\omega(e) = \phi(e)$ for every $e \in E(G)$.
			\end{cor}

		\subsection{Structure of the paper}

			To recap, during most of the rest of the paper we concentrate on the fractional version of Erd\H{o}s's question. More precisely, we are interested in minimising the quantity $\fpack(G)$, which is the largest size of a monochromatic fractional triangle packing in $G$, over red-blue colourings $G$ of $K_n$. 

			Our main result here is \Cref{thm:frac-main-almost-extremal}, which characterises colourings $G$ for which $\fpack(G) \le n(n-1)/4$. The proof of \Cref{thm:frac-main-almost-extremal} breaks down into three parts. \Cref{lem:computer} essentially tells us that we can focus on those $G$ that are either close to pentagon blow-ups or one of whose colour classes is close to bipartite. This lemma is proved by computer, and the algorithm is described in \Cref{sec:computer-search}. \Cref{lem:almost-bip-step} resolves the almost bipartite case, \Cref{lem:pentagon} resolves the pentagon blow-up case, and they are proved in \Cref{sec:almost-bip,sec:pentagon}, respectively. As mentioned at the end of \Cref{subsec:almost-extremal}, the proof of \Cref{thm:frac-main-almost-extremal} follows directly from \Cref{lem:computer,lem:almost-bip-step,lem:pentagon}.

			The next objective is to determine the minimum of $\fpack(G)$ over red-blue colourings of $K_n$. This is done in \Cref{thm:frac-main-extremal}, which is proved in \Cref{sec:frac-extremal}. The proof relies on the results described in the previous paragraph, as well as intermediate reults towards them. Together with the reduction of the original problem to its fractional version, due to Haxell and R\"odl \cite{haxell-rodl}, \Cref{thm:frac-main-extremal} implies the main result of this paper, namely \Cref{thm:main}, which confirms Erd\H{o}s's conjecture.

			Our final goal is to prove a stability version of our main result, \Cref{thm:main-stability}, which is what we do in \Cref{sec:stability}. To do so, we prove a `stability version' of our main fractional result (see \Cref{thm:frac-main-stability}), whose proof uses the main fractional result itself together with variants of arguments leading up to it. To deduce \Cref{thm:main-stability} from this `fractional stability result', following Tyomkyn \cite{tyomkyn}, we use a result of Alon, Shapira and Sudakov \cite{alon-shapira-sudakov}.

			We conclude the paper in \Cref{sec:conc} with some remarks and open problems.

	\section{Computer search} \label{sec:computer-search}

		In this section we describe the algorithm behind our computer search, which is used to prove \Cref{lem:computer}. The relevant certificates can be found \href{https://liveuclac-my.sharepoint.com/:f:/g/personal/ucahsle_ucl_ac_uk/EltPNwoeZx5GoXvNBgA3RDcBqGrsLjwjuZ1SW56SrlGK8A?e=5\%3ahICC1Q&at=9F}{\textcolor{blue}{here}}.
		
		Given a graph $G$ and a triangle packing $\omega$ in it, we define $\omega(G) = 3\sum_{T \in \T(G)} \omega(T)$, where $\T(G)$ is the collection of triangles in $G$.

		{\bf Initialisation.}\,
		Let $\LL_0$ be the set consisting of the empty graph (on $0$ vertices).

		{\bf Iteration.}\,
		Let $n \in \{0, \ldots, 21\}$, and let $\LL_n$ be a collection of red-blue colourings of $K_n$. 
		Initialise $\LL_{n+1} = \emptyset$.

		Given $H \in \LL_n$, let $u$ be a new vertex (so $u \notin V(H)$).
		We expose the edges from $u$ to $V(H)$ one by one, as follows.

		Start with the triple $(H, h_r, h_b)$, where $h_r$ and $h_b$ are maximum red and blue fractional triangle packings in $H$.

		Suppose that at some point we are given a triple $(F, f_r, f_b)$, where $f_r$ and $f_b$ are maximum red and blue fractional triangle packings in $F$. Proceed as follows.
		\begin{enumerate}
			\item \label{step:A}
				If $f_r(F) + f_b(F) > n(n+1)/4$, do nothing.
			\item \label{step:B}
				Otherwise, if $F$ is complete, proceed as follows. 
				\begin{enumerate}
					\item \label{step:Ba}
						If $n = 16$, check if $F$ is at most $1$-edge-flip away from a pentagon blow-up. If it is, take note of its blob sizes and whether it is a pentagon blow-up or $1$-away from a pentagon blow-up.
					\item \label{step:Bb}
						If $n = 21$, check if one of the colour classes in $F$ is $2$-close to bipartite. 
					\item \label{step:Bc}
						If $n \neq 16, 21$, or if $n = 16$ but $F$ is not close to a pentagon blow-up, or if $n = 21$ but neither colour in $F$ is close to bipartite, add $F$ to $\LL_{n+1}$, unless $\LL_{n+1}$ already contains a $2$-coloured graph isomorphic to $F$ or to the graph obtained by swapping the two colours in $F$.
				\end{enumerate}
			\item \label{step:C}
				Now suppose that $F$ is not complete, and that $f_r(F) + f_b(F) \le n(n+1)/4$. Pick a vertex $v$ such that $uv$ is not an edge in $F$, and proceed as follows.
				\begin{enumerate}
					\item \label{step:Ca}
						Form $F_r$ from $F$ by adding $uv$ as a red edge, and form $F_b$ by adding $uv$ as a blue edge.
					\item \label{step:Cb}
						Calculate a maximum fractional red triangle packing $f_r'$ in $F_r$, using a linear program with initial values set according to $f_r$. 
						Similarly, calculate a fractional blue triangle packing $f_b'$ in $F_b$.
					\item \label{step:Cc}
						Repeat the above procedure, with $(F_r, f_r', f_b)$ (note that $f_b$ is a maximum blue fractional triangle packing in $F_r$), and with $(F_b, f_r, f_b')$.
				\end{enumerate}
		\end{enumerate}

		{\bf Outcome.}\,
		The graphs found to be pentagon blow-ups when $n = 16$ in step \ref{step:Ba} have blobs sizes $x_1, \ldots, x_5$, where $(x_1, \ldots, x_5) \in \{(3,3,3,4,4), (2,3,4,4,4), (3,3,3,3,5)\}$; and the graphs found to be $1$-away from a pentagon blow-up have blob sizes $3, 3, 3, 4, 4$.

		At the end of the last iteration, when $n = 21$, $\LL_{22}$ remains empty.

		The proof of \Cref{lem:computer} readily follows.
		\begin{proof}[Proof of \Cref{lem:computer}]
			Let $G_0 \subseteq \ldots \subseteq G_{22}$ be as in the statement. If $G_{17}$ does not satisfy one of the first two items in the statement of \Cref{lem:computer} (namely, it is not close to a pentagon blow-up with blob sizes as in the outcome), then $G_{17} \in \LL_{17}$. If $G_{22}$ does not satisfy one of the last two items in the statement of the lemma, namely neither colour class in $G_{22}$ is $2$-close to bipartite, then $G_{22} \in \LL_{22}$, contrary to $\LL_{22}$ being empty.
		\end{proof}

		\subsection{Remarks}

			Here are some remarks regarding some technical aspects of our algorithm.

			\begin{enumerate} 
				\item
					In step \ref{step:Ba}, in order to determine if there is a $2$-coloured graph in $\LL_n$ which is isomorphic to $F$, we adapt an algorithm of McKay and Piperno \cite{mckay-piperno}.
				\item
					Throughout the process which explores the extensions of a given $H$, we use floats to store the weights of the triangles. 
					As such, the values of $f_b(F)$, $f_r(F)$, etc.\ are susceptible to rounding errors. 
					Thus, in step \ref{step:A}, before deciding to ignore $F$, we find a rational approximation of $f_r$ and $f_b$, using continuous fractions approximations (while ensuring that the weights are non-negative, and that no edge receives weight larger than $1$). We only ignore $F$ if these rational approximations of $f_r$ and $f_b$ still give a value larger than $n(n+1)/4$. 
					
					In practice, our program did not encounter such issues: whenever $f_r(F) + f_r(F)$ was found to be larger than $n(n+1)/4$, then the same held for the rational approximations. Nevertheless, to ensure correctness, this had to be checked. 

					In other words, the above rounding procedure ensures that there are no `false negatives'. The program is not, however, guaranteed to avoid `false positives'. Namely, some graphs $G \in \LL_n$ may theoretically not satisfy $\fpack(G) \le n(n-1)/4$, though this is unlikely in practice. This does not affect correctness.
				\item
					Similarly, the program that checks if a graph is close to bipartite may produce `false positives', but it is guaranteed not to produce `false negatives'. 
					In other words, every graph determined to be $2$-close to bipartite is indeed so (this is easy to verify given a suitable bipartition), but in theory our program could fail to find a suitable bipartition for a graph which is $2$-close to bipartite. In practice, as $\LL_{22}$ is empty at the end of the process, the algorithm for the almost bipartite case is successful on all the relevant graphs. 
				\item
					We note that our program that checks if a graph is $1$-close to a pentagon blow-up is always correct, namely it returns neither `false positives' nor `false negatives'. 
				\item
					The main element of the algorithm that allows us to improve on previous work (\cite{keevash-sudakov,tyomkyn}) is the fact that we expose the edges of the extension of $H$ one at a time, rather then exploring each of the extensions of $H$ into a red-blue colouring of a complete graph on $n+1$ vertices separately, as in \cite{keevash-sudakov,tyomkyn}.
				\item
					In step \ref{step:C} we choose the vertex $v$ following a simple greedy strategy. This speeds up the process considerably in comparison with a procedure that chooses $v$ according to a predetermined order of the vertices. In other words, graphs $F$ in step \ref{step:A} have typically fewer edges if $v$ is chosen judiciously.
				\item
					Using the values of $f_r$ as a basis for calculating $f_r'$ in step \ref{step:Cb}, introduces another important improvement by our algorithm. Indeed, intuitively, if two graphs differ by one edge, one could expect them to have similar maximum fractional triangle packings.
			\end{enumerate}

	\section{Almost bipartite} \label{sec:almost-bip}

		In this section we prove \Cref{lem:almost-bip-step}, which implies that if $G$ is a red-blue colouring of $K_{n+1}$ satisfying $\fpack(G) \le n(n+1)/4$, and, for some vertex $u$, $H = G \setminus \{u\}$ satisfies $\fpack(H) \le n(n-1)/4$ and $H_B$ is close to bipartite, then $G_B$ is close to bipartite. 
		The main ingredients in our proof are \Cref{thm:n-four} about fractional triangle packings in almost complete graphs, and \Cref{prop:min-four-applicable,prop:cover-k-triangles} (stated below) about red-blue colourings $H$ of $K_n$ with $\fpack(H) \le n(n-1)/4$, where $H_B$ is $(n/8)$-close to bipartite.

		The following proposition gives lower bounds on the size of the parts of the bipartition corresponding to $H_B$ being close to bipartite.

		\begin{prop} \label{prop:min-four-applicable}
			Let $n \ge 19$. Let $H$ be a red-blue colouring of $K_n$ satisfying $\fpack(H) \le n(n-1)/4$. Suppose that $\{X_1, X_2\}$ is a bipartition of $V(H)$ such that there are at most $n/8$ blue edges with both ends in either $X_1$ or $X_2$. Then the following two inequalities hold.
			\begin{enumerate} [label = \rm(\alph*)]
				\item \label{itm:X-lower-k-four}
					$|X_i| \ge k + 4$ for $i \in [2]$,
				\item \label{itm:X-lower-six}
					$|X_i| \ge 7$ for $i \in [2]$.
			\end{enumerate}
		\end{prop}

		Given a bipartition $\{X_1, X_2\}$ of $V(H)$, a \emph{cross triangle} is a triangle with at least one vertex in each of $X_1$ and $X_2$. The following proposition shows that there is a triangle packing in $H$, consisting of blue cross triangles, that covers the blue edges within $X_1$ and $X_2$. It will be useful later to find such a packing that avoids a given matching (which corresponds to blue triangles in $G$ containing $u$, where $H = G \setminus \{u\}$).

		\begin{prop} \label{prop:cover-k-triangles}
			Let $n \ge 22$. Let $H$ be a red-blue colouring of $K_n$ satisfying $\fpack(H) \le n(n-1)/4$. Suppose that $\{X_1, X_2\}$ is a bipartition of $V(H)$ such that there are at most $n/8$ blue edges with both ends in $X_1$ or in $X_2$.
			Let $M$ be a matching between $X_1$ and $X_2$. Then there is a triangle packing in $H \setminus M$, consisting of blue cross triangles, that covers all blue edges in $X_1$ or $X_2$. 
		\end{prop}

		We first prove \Cref{prop:min-four-applicable,prop:cover-k-triangles}, and then prove \Cref{thm:frac-main-almost-extremal}.

		\subsection{Proof of \Cref{prop:min-four-applicable}}

			\begin{proof}
				Write $|X_1| = n/2 + x$, so $|X_2| = n/2 - x$. Without loss of generality $|X_1| \ge |X_2|$, so it suffices to prove the required inequalities for $X_2$.

				First, we note that $|X_2| \ge 3$. Indeed, otherwise $H_R[X_2]$ is a graph on at least $n - 2$ vertices, with at most $n/8$ missing edges. By \Cref{thm:n-four}, as $n/8 \le n - 6$  (which holds for $n \ge 7$) and $n - 2 \ge 7$,  $H_R[X_1]$ has a red fractional triangle decomposition, whose size is at least $\binom{n-2}{2} - n/8$. As $\fpack(H) \le n(n-1)/4$, we have
				\begin{equation*} 
					n^2/4 - n/4 \ge \binom{n-2}{2} - n/8 = n^2/2 - 5n/2 + 3 - n/8 = n^2/2 - 21n/8 + 3.
				\end{equation*}
				This rearranges to $2n^2 - 19n + 24 \le 0$, a contradiction provided that $n \ge 9$. 

				Consider any fractional triangle decomposition of $X_1^{(2)} \cup X_2^{(2)}$ (as $|X_1| \ge |X_2| \ge 3$, such a decomposition exists), and remove all triangles that contain at least one blue edge. We end up with a red fractional triangle packing of size at least
				\begin{equation*} 
					\binom{n}{2} - |X_1| \cdot |X_2| - 3k 
					\ge \binom{n}{2} - (n/2 + x) \cdot (n/2 - x) - 3n/8 
					= n^2/4 - 7n/8 + x^2.
				\end{equation*}
				Using $\fpack(H) \le n(n-1)/4$ and rearranging, it follows that $x^2 \le 5n/8$.

				Now suppose that $n/2 - x < k + 4$, so $n/2 - x \le k + 3 \le n/8 + 3$. It follows that $x \ge 3n/8 - 3$. Note that $3n/8 - 3 \ge 1$, as $n \ge 11$, so $x^2 \ge (3n/8 - 3)^2$. Together with the inequality $x^2 \le 5n/8$, we find that 
				\begin{equation*} 
					5n/8 \ge (3n/8 - 3)^2 = 9n^2/64 - 9n/4 + 9.
				\end{equation*}
				Rearranging, this yields $9n^2 - 184n + 576 \le 0$, a contradiction for $n \ge 17$. This establishes \ref{itm:X-lower-k-four}.
				 
				Notice that if $k \ge 3$, \ref{itm:X-lower-six} follows from \ref{itm:X-lower-k-four}. We may thus assume that $k \le 2$ and that $|X_2| \le 6$. As before, by taking a fractional triangle decomposition in $X_1^{(2)} \cup X_2^{(2)}$ and removing triangles containing blue edges, we obtain a red fractional triangle packing of size at least 
				\begin{equation*} 
					\binom{n}{2} - |X_1| \cdot |X_2| - 6
					\ge n^2/2 - n/2 - 6(n-6) - 6 
					= n^2/2 - 6.5n + 30.
				\end{equation*}
				Using $\fpack(H) \le n(n-1)/4$ and rearranging, we find that $n^2 - 25n + 120 \le 0$, a contradiction for $n \ge 19$.
			\end{proof}

		\subsection{Proof of \Cref{prop:cover-k-triangles}}

			\begin{proof}
				Let $H' = H \setminus M$. Throughout this proof, we work with $H'$.

				Let $k_i$ be the number of blue edges in $X_i$, for $i \in [2]$; so $k = k_1 + k_2 \le n/8$.
				Let $\T_B$ be a largest triangle packing consisting of blue cross triangles. Let $m_i$ be the number of triangles in $\T_B$ with two vertices in $X_i$, for $i \in [2]$; so $m_i$ is the number of blue edges in $X_i$ that are covered by $\T_B$, and thus $m_i \le k_i$. Write $m = m_1 + m_2$. Our aim is to show that $m = k$, so suppose to the contrary that $m < k$.

				Let $\T_R$ be a largest \emph{fractional} triangle packing consisting of red cross triangles. If $|\T_B| + |\T_R| > n/2 - x - 4$ (where $|\T|$ denotes the number of triangles in $\T$, or the sum of weights of triangles), we reduce the weights of some triangles in $\T_R$ so that $|\T_R| = n/2 - x - 4 - |\T_B| = n/2 - x - 4 - m$ (note that this value is non-negative by \Cref{prop:min-four-applicable}~\ref{itm:X-lower-k-four}). Thus, by \Cref{cor:n-four-frac}, there is a fractional triangle decomposition in $H'_R[X_i] \setminus \T_R$; namely, the weighted graph obtained from $H'_R[X_i]$ by reducing the weight of each edge by its total weight in $\T_R$ has a fractional triangle decomposition. Combining these decompositions (for $i \in [2]$) with $\T_B$ and $\T_R$, we obtain a monochromatic fractional triangle packing in $H'$ of size at least
				\begin{equation*}
					\binom{n}{2} - (n/2 + x) \cdot (n/2 - x) - k + 3m + 2|\T_R|
					= n^2/4 - n/2 + x^2 -k + 3m + 2|\T_R|.
				\end{equation*}
				As $\fpack(H') \le \fpack(H) \le n(n-1)/4$, it follows that
				\begin{equation} \label{eq:packing}
					2|\T_R| - n/4 + 3m - k + x^2 \le 0.
				\end{equation}

				Suppose that $|\T_R| = n/2 - x - 4 - m$. Plugging this into \eqref{eq:packing}, we obtain
				\begin{equation*} 
					0 \ge 3n/4 + m - k + x^2 - 2x - 8
					\ge 5n/8 + (x - 1)^2 - 9
					\ge 5n/8 - 9,
				\end{equation*}
				a contradiction, as $n \ge 15$.
				So from now on we may assume that $\T_R$ is a largest fractional triangle packing consisting of cross red triangles.

				\begin{claim} \label{claim:both-sides-unsaturated}
					If $m_1 < k_1$ and $m_2 < k_2$, then $2|\T_R| \ge n - 2m - k - 10$.
				\end{claim}

				\begin{proof} 
					Let $T_i$ be the set of vertices in $X_i$ that appear in at least one triangle in $\T_B$, and let $a_i b_i$ be a blue edge in $X_i$ not covered by $\T_B$, for $i \in [2]$.  We note that $a_i$ and $b_i$ do not have a common blue neighbour in $X_{3-i} \setminus T_{3-i}$, as such a neighbour would give rise to a blue triangle that could be added to $\T_B$, contradicting the maximality of $\T_B$.
					For every blue edge in $X_i$ apart from $a_i b_i$, whose two ends are not in $T_i$, pick one of its vertices; let $S_i$ be the set of vertices chosen. So $|S_i| \le k_i - m_i - 1$ and $X_i \setminus (T_i \cup S_i)$ is a complete red graph.

					Let $A_i$ be the set of red neighbours of $a_{3-i}$ in $X_{i} \setminus (S_{i} \cup T_{i} \cup \{a_{i}, b_{i}\})$, and let $B_i$ be the set of red neighbours of $b_{3-i}$ in $X_{i} \setminus (A_i \cup S_{i} \cup T_{i} \cup \{a_{i}, b_{i}\})$.
					As $a_{3-i}$ and $b_{3-i}$ do not have common blue neighbours in $X_{i} \setminus T_{i}$, and they each have at most one non-neighbour in $X_i$, we have $|A_i| + |B_i| \ge |X_i| - |T_i| - |S_i| - 4$.

					As $A_i$ and $B_i$ span complete red graphs, they span red matchings that cover at least $|A_i| - 1$ and $|B_i| - 1$ vertices, respectively. Consider the red triangle packing formed by taking these matchings in $A_i$ and $B_i$ and attaching their edges to $a_{3-i}$ and $b_{3-i}$, respectively, for $i \in [2]$. This is a red triangle packing $\T$ with
					\begin{align*}
						2|\T| \ge |A_1| + |B_1| + |A_2| + |B_2| - 4 
						& \ge |X_1| + |X_2| - |T_1| - |T_2| - |S_1| - |S_2| - 12 \\
						& \ge n - 3m - (k - m - 2) - 12 \\
						& = n - 2m - k - 10.
					\end{align*}
					The claim follows by maximality of $\T_R$.
				\end{proof}

				\begin{claim} \label{claim:one-side-unsaturated}
					If $m_i < k_i$ for exactly one $i \in [2]$, then $2|\T_R| \ge n/2 - x - 3m - 2 - \one\{k \le 2\}$.
				\end{claim}

				\begin{proof} 
					Let $\sigma \in [2]$ be such that $m_{\sigma} = k_{\sigma}$; write $\tau = 3 - \sigma$.
					Let $T$ be the set of vertices in $X_{\sigma}$ that appear in at least one of the triangles in $\T_B$; so $|T| \le 2m$.

					Let $ab$ be a blue edge in $X_{\tau}$ that is not covered by $\T_B$.
					Let $A$ be the set of red neighbours of $a$ in $X_{\sigma} \setminus T$, and let $B$ be the set of red neighbours of $b$ in $X_{\sigma} \setminus (T_{\sigma} \cup A)$. Then 
					\begin{equation} \label{eq:AB}
						|A| + |B| \ge |X_{\sigma}| - |T| - 2 \ge n/2 - x - 2m - 2.
					\end{equation}
					Note that if $|A| \neq 1$, then $A$ has a perfect (red) fractional matching (e.g.\ take the single edge in $A$ if $|A| = 2$, and otherwise give weight $1/2$ to each of the edges of a Hamilton cycle); a similar statement holds for $B$. Thus, attaching $a$ and $b$ to the edges of a perfect fractional matching in $A$ and $B$, respectively, we obtain a fractional triangle packing $\T$ consisting of red cross triangles, with 
					\begin{align*} 
						2|\T|
						& \ge |A| + |B| - \one\{|A| = 1\} - \one\{|B| = 1\} \\
						& \ge n/2 - x - 2m - 2 - \one\{|A| = 1\} - \one\{|B| = 1\} \\
						& \ge n/2 - x - 2m - 4.
					\end{align*}
					If $m \ge 2$, this is at least $n/2 - x - 3m - 2$, as required.
					If $m = 1$, then by \eqref{eq:AB} and \Cref{prop:min-four-applicable} $|A| + |B| \ge 3$, so at most one of $A$ and $B$ has size $1$. It follows that $2|\T| \ge n/2 - x - 2m - 3 = n/2 - x - 3m - 2$, as required.
					If $m = 0$, then again at most one of $A$ and $B$ has size $1$. We are done if both $A$ and $B$ do not have size $1$; or if $k \le 2$; or if $|A| + |B| > n/2 - x - 2$. So we may assume that $|B| = 1$, $|A| = n/2 - x - 3 \ge 4$, and $k \ge 3$. In fact, we may assume that for every blue edge $a'b'$ in $X_{\tau}$, if $A'$ and $B'$ are defined as above, then (after possibly swapping the roles of $a'$ and $b'$) $|B'| = 1$ and $|A'| = n/2 - x - 3 \ge 4$.

					Suppose that there are two vertices $a$ and $a'$ in $X_{\tau}$ with $n/2 - x - 3$ red neighbours in $X_{\sigma}$, and let $A$ and $A'$ be their red neighbourhoods in $X_{\sigma}$. Form a fractional triangle packing $\T'$, consisting of red cross triangles, by attaching the edges of Hamilton cycles in $A$ and in $A'$ to $a$ and $a'$, respectively, giving each triangle weight $1/2$ (note that this is indeed a fractional triangle packing). Then $2|\T'| \ge |A| + |A'| \ge 2(n/2 - x - 3) \ge n/2 - x - 3m + 1$, say, and the claim readily follows. We may thus assume that there is only one vertex $a \in X_{\tau}$ with $n/2 - x - 3$ red neighbours in $X_{\sigma}$.

					As $X_{\tau}$ contains $k \ge 3$ blue edges, the above assumptions imply that there are three vertices $b_1, b_2, b_3 \in X_{\tau}$ such that $ab_i$ is blue. These assumptions also imply that there is a vertex $c \in X_{\sigma}$ such that $b_1c$ is a red edge and $ac$ is a blue edge. Without loss of generality, $b_2c$ is an edge (as every vertex has at most one non-neighbour). If $b_2c$ is red, we can add the triangle $b_1b_2c$ to $\T$ (with weight $1$), thus obtaining a fractional triangle packing of red cross triangles of the required size. Otherwise, $b_2c$ is blue, but then $ab_2c$ is a blue triangle, contradicting $m = 0$. 
				\end{proof}
				
				By \Cref{claim:both-sides-unsaturated,claim:one-side-unsaturated}, it suffices to consider two cases: $2|\T_R| \ge n - 2m - k - 10$; and $2|\T_R| \ge n/2 - x - 3m - 2 - \one\{k \le 2\}$.
				In the first case, we obtain the following inequality, using \eqref{eq:packing}
				\begin{equation*} 
					0 \ge 3n/4 + m - 2k + x^2 - 10
					\ge n/2 - 10,
				\end{equation*}
				a contradiction to $n \ge 21$.
				In the second case, we obtain the following, using \eqref{eq:packing}
				\begin{align*} 
					0 
					& \ge n/4 - k + x^2 - x - 2 - \one\{k \le 2\} \\
					& = n/4 - k + (x - 1/2)^2 - 2.25  - \one\{k \le 2\} \\
					& \ge n/4 - k - 2.25 - \one\{k \le 2\}.
				\end{align*}
				If $k \ge 3$, this gives $0 \ge n/8 - 2.25$, a contradiction if $n \ge 19$; and if $k \le 2$, this gives $0 \ge n/4 - 5.25$, a contradiction if $n \ge 22$.
			\end{proof}

		\subsection{Proof of \Cref{lem:almost-bip-step}}

			\begin{proof}
				Let $\{X_1, X_2\}$ be a bipartition of $V(H)$ such that the number of blue edges with both ends in $X_1$ or in $X_2$ is $k \le n/8$; denote by $k_i$ the number of blue edges in $X_i$ for $i \in [2]$. Write $|X_1| = n/2 + x$, so $|X_2| = n/2 - x$, and we assume that $x \ge 0$.
				Let $R_i$ and $B_i$ be the red and blue neighbourhoods of $u$ in $X_i$, respectively, for $i \in [2]$.

				\begin{claim} \label{claim:full-matching}
					There is a blue matching between $B_1$ and $B_2$ that saturates the smaller of the two sets.
				\end{claim}

				\begin{proof} 
					Let $M$ be a matching of maximum size between either $B_1$ and $B_2$, denote its size by $m$, and suppose that it does not saturate $B_1$ nor $B_2$. Write $A_i = B_i \setminus V(M)$; then $|A_i| \ge 1$ for $i \in [2]$. Let $a_i$ be any vertex in $A_i$.
					Let $\T_{R, 1}$ be a red triangle packing formed by attaching $a_{2}$ to a maximum matching in $A_1 \setminus \{a_1\}$, attaching $u$ to a maximum matching in $R_1$, and removing the triangles that contain a blue edge. Similarly, let $\T_{R,2}'$ be a red triangle packing formed by attaching $a_1$ to a maximum matching in $A_2$, attaching $u$ to a maximum matching in $A_2$, and removing triangles containing a blue edge. 
					Note that the triangles in $\T_{R,1}$ and $\T_{R,2}'$ are edge-disjoint, and
					\begin{align*} 
						& 2|\T_{R, 1}| \ge |R_1| - 1 + |A_1| - 2 - 2k_1 = |X_1| - m - 3 - 2k_1 \\
						& 2|\T_{R, 2}'| \ge |R_2| - 1 + |A_2| - 1 - 2k_2 = |X_2| - m - 2 - 2k_2.						
					\end{align*}
					We note that $|\T_{R, 1}| \le |X_1| - k_1 - 4$. Indeed, otherwise, as $|\T_{R, 1}| \le (|X_1 - 1)/2$, we get $(|X_1| - 1)/2 \ge |X_1| - k_1 - 3$, which rearranges to
					\begin{equation*} 
						0 \ge |X_1|/2 - k_1 - 5/2
						\ge n/4 - n/8 - 5/2 
						= n/8 - 5/2,
					\end{equation*}
					a contradiction to $n \ge 21$.
					Let $\T_{R,2}$ be a subpacking of $\T_{R,2}'$ of size $|\T_{R, 2}| = \min\{|X_2| - k_2 - 4, |\T_{R,2}'|\}$ (note that $|X_2| \ge k_2 + 4$ by \Cref{prop:min-four-applicable}), and let $\T_R = \T_{R, 1} \cup \T_{R, 2}$.

					Let $\T_B'$ be a blue triangle packing in $H \setminus M$ that consists of blue cross triangles, and that covers all blue edges in $X_1$ and in $X_2$; such a packing exists by \Cref{prop:cover-k-triangles}. Let $\T_B''$ be the blue triangle packing in $G$ obtained by attaching $u$ to the edges of $M$. 

					By choice of $\T_R$ and by \Cref{thm:n-four}, there is a fractional triangle decomposition in $H_R[X_i] \setminus \T_R$.
					It follows that there exists a monochromatic triangle packing in $G$ of size at least
					\begin{align*} 
						\binom{n}{2} - (n/2 + x) \cdot (n/2 - x) + 2|\T_B'| + 3|\T_B''| + 2|\T_R|
						= n^2/4 - n/2 + x^2 + 2k + 3m + 2|\T_R|.
					\end{align*}
					As $\fpack(G) \le n(n+1)/4$, we obtain 
					\begin{equation} \label{eq:packing-two}
						2|\T_R| - 3n/4 + 3m + 2k + x^2 \le 0.
					\end{equation}
					Recall that either $2|\T_R| = 2|\T_{R,1}'| + 2|\T_{R,2}'| \ge n - 2m - 5 - 2k$; or $2|\T_R| = 2|\T_{R,1}| + 2(|X_2| - k_2 - 4) \ge 3n/2 - m - 2k - x - 11$.
					If the former holds, then, by \eqref{eq:packing-two},
					\begin{align*} 
						0 \ge n/4 + m + x^2 - 5 \ge n/4 - 5, 
					\end{align*}
					a contradiction to $n \ge 21$.
					If the latter holds, then, again by \eqref{eq:packing-two},
					\begin{equation*} 
						0 \ge 3n/4 + 2m + x^2 - x - 11 = 3n/4 + 2m + (x-1/2)^2 - 11.25 \ge 3n/4 - 11.25,
					\end{equation*}
					a contradiction to $n \ge 16$.
				\end{proof}
				 
				Let $m = \min\{|B_1|, |B_2|\}$. 
				Let $\sigma \in [2]$ be such that $|B_{\sigma}| = m$. Define $X_{\sigma}' = X_{\sigma} \cup \{u\}$ and $X_{3 - \sigma}' = X_{3 - \sigma}$.
				In order to complete the proof, we need to show that the number of blue edges with both ends in $X_1'$ or in $X_2'$ is at most $(n+1)/8$; in other words, our task is to show that $m + k \le (n+1)/8$. 
				
				Let $M$ be a blue matching of size $m$ between $B_1$ and $B_2$; by \Cref{claim:full-matching}, such a matching exists.
				Let $\T_B'$ be a triangle packing in $H \setminus M$, that consists of blue cross triangles, and that covers all blue edges in $X_1$ or $X_2$; by \Cref{prop:cover-k-triangles}, such a packing exists. Let $\T_B''$ be the blue triangle packing obtained by attaching $u$ to each of the edges of $M$. Note that the triangles in $\T_B'$ and $\T_B''$ are edge-disjoint.

				\begin{claim} \label{claim:m-small}
					$m + k \le n/2 - x - 4$.
				\end{claim}

				\begin{proof} 
					Suppose not. Then, as $n/2 - x$ is an integer, we have 
					\begin{equation} \label{eq:m}
						m + k \ge n/2 - x - 3.
					\end{equation}
					By \Cref{thm:n-four} and \Cref{prop:min-four-applicable} there is a fractional triangle decomposition in $H_R[X_i]$ for $i \in [2]$. This packing, together with $\T_B'$ and $\T_B''$, forms a monochromatic fractional triangle packing in $G$ of size at least
					\begin{equation*} 
						\binom{n}{2} - (n/2 + x) \cdot (n/2 - x) + 2|\T_B'| + 3|\T_B''|
						= n^2/4 - n/2 + x^2 + 2k + 3m.
					\end{equation*}
					As $\fpack(G) \le n(n+1)/4$, we have
					\begin{align*} 
						0 \ge -3n/4 + x^2 + 2k + 3m
						& = -3n/4 + x^2 + 3(k + m) - k \\
						& \ge 3n/4 + x^2 - 3x - 9 - k \\
						& \ge 5n/8 + (x - 3/2)^2 - 11.25 \\
						& \ge 5n/8 - 11.25,
					\end{align*}
					using \eqref{eq:m} for the second inequality. This is a contradiction to $n \ge 20$.
				\end{proof}

				By \Cref{thm:n-four}, \Cref{prop:min-four-applicable} and \Cref{claim:m-small} there is a fractional triangle decomposition in $G_R[X_i']$ for $i \in [2]$. Together with $T_B'$ and $T_B''$ this forms a triangle packing in $G$ of size at least
				\begin{align*} 
					\binom{n+1}{2} - |X_1'| \cdot |X_2'| + 2(k + m)
					& \ge n(n+1)/2 - (n+1)^2/4 + 2(k + m) \\
					& \ge (n+1)(n-1)/4 + 2(k+m).
				\end{align*}
				As $\fpack(G) \le n(n+1)/4$, it follows that $k+m \le (n+1)/8$, as required.
			\end{proof}

	\section{Minimising the size of a fractional triangle packing} \label{sec:frac-extremal}

		In this section we prove \Cref{thm:frac-main-extremal}, which determines the minimum of $\fpack(G)$ among red-blue colourings $G$ of $K_n$, for $n \ge 26$, and characterises the minimisers.

		\begin{proof}[Proof of \Cref{thm:frac-main-extremal}]
			Let $G$ be a red-blue colouring of $K_n$ where $G_B$ is obtained by removing a matching from the complete bipartite graph $K_{\floor{n/2}, \ceil{n/2}}$. Then $G_B$ is a triangle-free. The (red) edges within the two parts can all be covered by a red fractional triangle packing (as each of the two parts forms a complete red graph on more than two vertices), and the red edges between the parts are not contained in any triangle. It follows that 
			\begin{equation*} 
				\fpack(G) 
				= \binom{n}{2} - \floor{n/2} \cdot \ceil{n/2} 
				= n^2/2 - n/2 - \ceil{\frac{n^2 - 1}{4}}
				= \floor{\frac{(n-1)^2}{4}}.
			\end{equation*}

			Now suppose that $G$ is a red-blue colouring of $K_n$ with $\fpack(G) \le \floor{(n-1)^2/4}$. By \Cref{thm:frac-main-almost-extremal}, without loss of generality, there is a bipartition $\{X_1, X_2\}$ of the vertices of $G$ such that there are $k \le n/8$ blue edges within the parts $X_1$ and $X_2$. 
			Without loss of generality, $|X_1| \ge |X_2|$. Write $|X_1| = \ceil{n/2} + x$, so $|X_1| = \floor{n/2} - x$.

			By \Cref{prop:cover-k-triangles} there is a triangle packing $\T_B$ that consists of blue cross triangles and covers all blue edges within the parts $X_1$ and $X_2$. By \Cref{prop:min-four-applicable} and \Cref{thm:n-four}, $G_R[X_1]$ and $G_R[X_2]$ have fractional triangle decompositions. Putting these together with $\T_B$, we obtain a monochromatic fractional triangle packing of size at least
			\begin{align*} 
				\binom{n}{2} - |X_1| \cdot |X_2| + 2|\T_B|
				& = \binom{n}{2} - \ceil{n/2} \cdot \floor{n/2} + x(\ceil{n/2} - \floor{n/2}) + x^2 + 2k \\
				& \ge \floor{\frac{(n-1)^2}{4}} + x^2 + 2k.
			\end{align*}
			As $\fpack(G) \le \floor{(n-1)^2/4}$, it follows that $x = k = 0$, so $|X_1| = \ceil{n/2}$, $|X_2| = \floor{n/2}$, and $G_R[X_1]$ and $G_R[X_2]$ are complete. Now consider $G_R[X_1, X_2]$. Suppose that there is a vertex $u$ with degree at least $2$ in this graph. Let $i \in [2]$, $u \in X_i$ and $v, w \in X_{3-i}$ be such that $v$ and $w$ are red neighbours of $u$. Again by \Cref{prop:min-four-applicable} and \Cref{thm:n-four}, $G_R[X_i]$ and $G_R[X_{3-i}] \setminus \{vw\}$ both have triangle decompositions. It follows that 
			\begin{equation*} 
				\fpack(G) \ge \binom{n}{2} - \ceil{n/2} \cdot \floor{n/2} + 2 \ge \floor{\frac{(n-1)^2}{4}} + 2,
			\end{equation*}
			a contradiction. So $G_R[X_1, X_2]$ has maximum degree at most $1$, i.e.\ $G_B$ is $K_{\ceil{n/2}, \floor{n/2}}$ minus a matching, as required.
		\end{proof}
	\section{Stability} \label{sec:stability}
		
		In this section we prove our second main theorem, \Cref{thm:main-stability}, which is a stability result regarding monochromatic triangle packings in $2$-coloured complete graphs. As mentioned in \Cref{sec:overview}, we will deduce it from \Cref{thm:frac-main-stability}, which we will prove later in this section.

		For an $n$-vertex graph $G$, let $E_{\bip}(G)$ be the least $\delta$ such that $G$ is $\delta n^2$-close to bipartite, i.e.\ the least $\delta$ such that $G$ can be made bipartite by the removal of at most $\delta n^2$ edges. We use the following special case of a result by Alon, Shapira and Sudakov \cite{alon-shapira-sudakov}.

		\begin{thm}[A special case of Theorem 1.2 in \cite{alon-shapira-sudakov}] \label{thm:a-s-s}
			For every $\mu > 0$ there exists $d = d(\mu)$ with the following property: let $G$ be a graph, let $D$ be a subset of $V(G)$ of size $d$, chosen uniformly at random among such sets. Then
			\begin{equation*} 
				\Pr\Big[ \left|E_{\bip}(G) - E_{\bip}(G[D])\right| > \mu \Big] < \mu.
			\end{equation*}
		\end{thm}
		It is not hard to convince oneself that $d(\mu)$ tends to infinity as $\mu$ goes to $0$.

		\begin{proof}[Proof of \Cref{thm:main-stability} using \Cref{thm:frac-main-stability}]
			Let $\eps > 0$, let $\mu > 0$ be sufficiently small, and let $d = d(\mu)$ be as given by \Cref{thm:a-s-s}.
			Let $G$ be a red-blue colouring of $K_n$ where both $G_B$ and $G_R$ are $\eps n^2$-far from bipartite, i.e.\ $E_{\bip}(G_R), E_{\bip}(G_B) \ge \eps$.
			By \Cref{thm:a-s-s}, for all but at most $2\mu \binom{n}{d}$ sets $D$ of $d$ vertices, both $G_R[D]$ and $G_B[D]$ are $(\eps - \mu)d^2$-far from bipartite. Assuming that $\mu$ is sufficiently small, and thus $d$ is sufficiently large, $(\eps - \mu)d^2 \ge (\eps/2)d^2 \ge (1/8 + \eta)d$, where $\eta$ is as in \Cref{thm:frac-main-stability}. By \Cref{thm:frac-main-stability}, $\fpack(G[D]) \ge d(d-1)/4 + 2\eta d$ for all but at most $2\eta \binom{n}{d}$ sets of $d$ vertices $D$. 
			Recall that $\fpack(G[D]) \ge d(d-2)/4 = d(d-1)/4 - d/4$ for all sets of $d$ vertices $D$, assuming $d \ge 22$, by \Cref{thm:frac-main-extremal}. 
			By iterating \Cref{obs:averaging}, we have
			\begin{align*} 
				\fpack(G) 
				& \ge \frac{1}{\binom{n-2}{d-2}} \sum_{D \subseteq V(G):\, |D| = d} \fpack(G[D]) \\
				& \ge \frac{1}{\binom{n-2}{d-2}} \cdot \binom{n}{d} \Big( (1 - 2\mu) \cdot \big(d(d-1)/4 + 2\eta d\big) + 2\mu \cdot \big(d(d-1)/4 - d/4\big)\Big) \\
				& \ge \frac{n(n-1)}{d(d-1)} \cdot \Big(d(d-1)/4 + d \cdot \big( (1 - 2\mu) \cdot 2\eta - \mu/2 \big)\Big) \\
				& \ge \frac{n(n-1)}{d(d-1)} \cdot \Big(d(d-1)/4 + \eta d \Big) \\
				& = \frac{n(n-1)}{4} + \eta \cdot \frac{n(n-1)}{d-1} \\
				& \ge \left(\frac{1}{4} + \frac{\eta}{2d}\right)n^2,
			\end{align*}
			where the inequalities hold for sufficiently small $\mu$ and sufficiently large $n$. \Cref{thm:main-stability} follows by taking $\delta = \eta / (2d)$.
		\end{proof}

		\subsection{Overview of the proof of \Cref{thm:frac-main-stability}}

			The following lemma, which is a strengthening of \Cref{lem:almost-bip-step} (for large $n$), is the main ingredient in the proof of \Cref{thm:frac-main-stability}.

			\begin{lem} \label{lem:bip-stability}
				The following holds for sufficiently large $n$.
				Let $G$ be a red-blue colouring of $K_{n+1}$, let $u \in V(G)$, and set $H = G \setminus \{u\}$.
				Suppose that $H_B$ is $(1/8 + 1/200)n$-close to bipartite and that $\fpack(H) \le n(n-1)/4 + n/200$.
				Then for any $\ell \le (1/8 + 1/200)(n+1)$, either $\fpack(G) \ge (n^2 - 1)/4 + 2\ell$, or $G_B$ is $\ell$-close to bipartite.
			\end{lem}

			\begin{proof}[Proof of \Cref{thm:frac-main-stability} using \Cref{lem:bip-stability}]
				Let $n_0$ be large enough so that \Cref{lem:bip-stability} is applicable for $n \ge n_0$.
				We prove the theorem by induction on $n \ge n_0$. For the base case $n = n_0$, recall that by \Cref{thm:frac-main-almost-extremal}, either one of $G_B$ and $G_R$ is $(n/8)$-close to bipartite, or $\fpack(G) > n(n-1)/4$. Graphs satisfying the former condition immediately satisfy the requirements of the theorem, and by taking $\eta \in (0, 1/200)$ to be sufficiently small, we can guarantee that all graphs $G$ satisfying the latter condition also satisfy $\fpack(G) \ge n(n-1)/4 + \eta n$. 

				Now take $n \ge n_0$, and suppose that the statement of \Cref{thm:frac-main-stability} holds for $n$. Let $G$ be a red-blue colouring of $K_{n+1}$. Then, by induction, for every vertex $u$ in $G$, either $\fpack(G \setminus \{u\}) \ge n(n-1)/4 + 2\eta n$, or one of $G_B \setminus \{u\}$ and $G_R \setminus \{u\}$ is $(1/8 + \eta)n$-close to bipartite. If the former holds for all vertices $u$, then 
				\begin{align*} 
					\fpack(G) 
					& \ge \frac{1}{n-1} \sum_{u \in V(G)} \fpack\left(G \setminus \{u\}\right) \\
					& \ge \frac{1}{n-1} \cdot (n+1) \cdot \left(\frac{n(n-1)}{4} + 2\eta n\right) \\
					& \ge \frac{n(n+1)}{4} + 2\eta(n+1),
				\end{align*}
				as required.
				So we may assume, without loss of generality, that for some vertex $u$, the colouring $H = G \setminus \{u\}$ satisfies $\fpack(H) \le n(n-1)/4 + 2\eta n$ and $H_B$ is $(1/8 + \eta)n$-close to bipartite.
				By \Cref{lem:bip-stability} (with $\ell = (1/8 + \eta)(n+1)$), either $\fpack(G) \le (n^2 - 1)/4 + (n+1)/4 + 2\eta(n+1) = n(n+1)/4 + 2\eta(n+1)$, or $G_B$ is $(1/8 + \eta)n$-close to bipartite.
			\end{proof}

			In order to prove \Cref{lem:bip-stability}, we make use of the following two propositions, which are variants of \Cref{prop:min-four-applicable,prop:cover-k-triangles}. 

			\begin{prop} \label{prop:bip-stability-x}
				Let $H$ be a red-blue colouring of $K_n$ satisfying $\fpack(H) \le n^2/4$. Let $\{X_1, X_2\}$ be a partition of the vertices of $H$ with at most $(n-2)/6$ blue edges with both ends in either $X_1$ or $X_2$. Then $n/2 - \sqrt{n} \le |X_i| \le n/2 + \sqrt{n}$.
			\end{prop}

			\begin{prop} \label{prop:bip-stability-blue-neighs}
				The following holds for sufficiently large $n$.
				Let $H$ be a red-blue colouring of $K_n$. Let $\{X_1, X_2\}$ be a bipartition of the vertices of $H$ such that there are at most $n/8 + n/200$ blue edges with both ends in either $X_1$ or $X_2$, and suppose that $\fpack(H) \le n(n-1)/4 + n/200$. Then 
				\begin{enumerate}[label = \rm(\alph*)]
					\item \label{itm:stability-few-red-neighs}
						every vertex in $X_i$ has at most $n/20$ red neighbours in $X_{3-i}$, for $i \in [2]$,
					\item \label{itm:stability-perfect-blue-packing}
						given any matching $M$ in $H[X_1, X_2]$, there is a triangle packing that consists of blue cross triangles avoiding the edges of $M$, and that covers all blue edges in either $X_1$ or $X_2$.
				\end{enumerate}
			\end{prop}

			We first prove the two propositions, and then deduce \Cref{lem:bip-stability} from them.

		\subsection{Proofs of \Cref{prop:bip-stability-x,prop:bip-stability-blue-neighs}}

			\begin{proof}[Proof of \Cref{prop:bip-stability-x}]
				Let $k$ be the number of blue edges with both ends in either $X_1$ or $X_2$. Consider a fractional triangle decomposition in $G[X_1]$ and in $G[X_2]$; note that such decompositions exist, unless one of $X_1$ or $X_2$ has size exactly $2$, in which case we leave one edge uncovered. Now remove all triangles that contain a blue edge. This results in a red fractional triangle packing of size at least the following, where $|X_1| = n/2 + x$ (and $|X_2| = n/2 - x$),
				\begin{equation*} 
					\binom{n}{2} - (n/2 + x) \cdot (n/2 - x) - 1 - 3k
					= n^2/4 - n/2 + x^2 - 1 - 3k.
				\end{equation*}
				As $\fpack(G) \le n^2/4$, it follows that $x^2 \le n/2 + 1 + 3k \le n$, as required.
			\end{proof}

			\begin{proof}[Proof of \Cref{prop:bip-stability-blue-neighs}]
				Let $k$ be the number of blue edges in $X_1$ or $X_2$.

				\begin{claim} \label{claim:common-blue-neighs}
					Every two vertices in $X_i$ have at least $(1/8 - 1/100)n - 3\sqrt{n}$ common blue neighbours in $X_{3-i}$, for $i \in [2]$.
				\end{claim}
				\begin{proof}
					We prove the statement for $i = 1$, the claim would follow by symmetry.
				
					Let $x_1, x_2 \in X_1$, and let $\ell$ be the number of common blue neighbours of $x_1$ and $x_2$. Let $R_1$ be the red neighbourhood of $x_1$ in $X_2$, and let $R_2$ be the red neighbourhood of $x_2$ in $X_2 \setminus R_1$; so $|R_1| + |R_2| \ge |X_2| - \ell$. Let $M_i$ be a maximum matching in $H_R[R_i]$. As the vertices uncovered by $M_i$ form an independent set, and there are at most $k$ missing edges, $M_i$ covers at least $|R_i| - \sqrt{2k} - 1$ vertices. 
					Let $\T_R$ be the red triangle packing obtained by attaching $x_i$ to the edges of $M_i$, for $i \in [2]$. 
					Note that $H_R[X_1]$ and $H_R[X_2] \setminus (M_1 \cup M_2)$ have fractional triangle decompositions; this follows from \Cref{thm:n-four}, as the number of edges missing from either graph is at most $k + |X_2|/2 \le n/2 - \sqrt{n} - 4 \le |X_1|, |X_2|$, using \Cref{prop:bip-stability-x}. 
					We conclude that $H$ has a red fractional triangle packing of size at least
					\begin{align*} 
						\binom{n}{2} - |X_1| \cdot |X_2| - k + 2(|M_1| + |M_2|)
						& \ge n^2/4 - n/2 - k + |R_1| + |R_2| - 2\sqrt{2k} + 1 \\
						& \ge n^2/4 - n/2 - k + n/2 - \sqrt{n} - \ell - 2\sqrt{n} \\
						& = n^2/4 - 3\sqrt{n} - k - \ell.
					\end{align*}
					As $\fpack(H) \le n(n-1)/4 + n/200$, it follows that $k + \ell \ge n/4 - n/200 - 3\sqrt{n}$, implying that $\ell \ge n/8 - 2n/200 - 3\sqrt{n}$, as required.
				\end{proof}

				Write $\ell = (1/8 - 1/100)n - 3\sqrt{n}$.
				Let $k_i$ be the number of blue edges in $X_i$, for $i \in [2]$, and let $m_i = \min\{k_i, \ell - 2\}$.
				Let $e_1, \ldots, e_{m_1}$ be distinct blue edges in $X_1$, and let $f_1, \ldots, f_{m_2}$ be distinct blue edges in $X_2$. Let $x_1, \ldots, x_{m_1}$ be distinct vertices in $X_2$ such that $x_i$ is a common blue neighbour of the two ends of $e_i$, for $i \in [m_1]$. Note that such vertices can be chosen greedily by the above claim. Similarly, let $y_1, \ldots, y_{m_2}$ be distinct vertices in $X_1$, such that $y_i$ is a common blue neighbour of $f_i$, and the edges between $y_i$ and the ends of $f_i$ do not contain any edge between $x_j$ and $e_j$ for $j \in [m_1]$ and $i \in [m_2]$. Again, such vertices can be chosen greedily due to \Cref{claim:common-blue-neighs}, noting that when $y_i$ is chosen, there are at most $i+1 \le \ell - 1$ `forbidden neighbours'; namely, $y_i$ needs to be distinct from $y_1, \ldots, y_{i-1}$ and from the at most two neighbours in $X_1$ of the two ends of $f_i$ following edges between $x_j$ and $e_j$, for $j \in [m_1]$. Let $\T_B$ be the blue triangle packing obtained by attaching $x_i$ to $e_i$ and $y_j$ to $f_j$ for $i \in [m_1]$ and $j \in [m_2]$. 

				Let $u \in X_1$, and let $r$ be the number of red neighbours of $u$ in $X_2$. We will show that $r \le n/20$, which, by symmetry, would establish \ref{itm:stability-few-red-neighs}. Let $R$ be the red neighbourhood of $u$ in $X_2$, so $|R| = r$. As in the proof of \Cref{claim:common-blue-neighs}, $H_R[R]$ contains a matching $M$ that covers at least $|R| - \sqrt{2k} - 1$ vertices. Let $\T_R$ be the red triangle packing obtained by attaching $u$ to each of the edges of $M$.

				By noting that $H_R[X_i] \setminus M$ and $H_R[X_{3-i}]$ have fractional triangle decomposition, and by considering the packings $\T_B$ and $\T_R$, it follows that 
				\begin{align*} 
					\fpack(H) 
					& \ge n^2/4 - n/2 - k + 3\min\{k, \ell - 2\} + r - \sqrt{2k} - 1 \\
					& \ge n^2/4 - n/2 - n/8 - n/200 + 3(n/8 - n/100 - 3\sqrt{n} - 2) + r + \sqrt{n} \\
					& = n^2/4 - n/4 - 7n/200 - 11\sqrt{n} + r.
				\end{align*}
				As $\fpack(H) \le n(n-1)/4 + n/200$, it follows that $r \le 8n/200 + 11\sqrt{n} \le n/20$, as required for \ref{itm:stability-few-red-neighs}.

				It remains to prove \ref{itm:stability-perfect-blue-packing}. Let $M$ be a matching in $H[X_1, X_2]$. Let $e_1, \ldots, e_{k_1}$ be the blue edges in $X_1$, and let $f_1, \ldots, f_{k_2}$ be the blue edges in $X_2$. As above, we find distinct vertices $x_1, \ldots, x_{k_1} \in X_2$ such that $x_i$ is a common blue neighbour of the ends of $e_i$, and neither of the edges between $x_i$ and $e_i$ are in $M$, for $i \in [2]$. As the number of `forbidden' vertices at any point is at most $k_1 \le k$ and the number of common neighbours of the ends of $e_i$ is at least $|X_2| - 2n/10 \ge k + 1$, such vertices $x_i$ exist. Similarly, we may find distinct vertices $y_1, \ldots, y_{k_2} \in X_1$ such that $y_i$ is a common blue neighbour of the ends of $f_i$, and neither of the edges between $y_i$ and $f_i$ are in $M$, or are an edge between $x_j$ and $e_j$ for $j \in [k_1]$ and $i \in [k_2]$. Attaching $x_i$ to $e_i$ for every $i \in [k_1]$ and attaching $y_j$ to $f_j$ for every $j \in [k_2]$, gives a triangle packing as required for \ref{itm:stability-perfect-blue-packing}.
			\end{proof}

		\subsection{Proof of \Cref{lem:bip-stability}}

			It remains to prove \Cref{lem:bip-stability}. This is now a fairly straightforward task.

			\begin{proof}[Proof of \Cref{lem:bip-stability}]
				Let $\{X_1, X_2\}$ be a bipartition of $H$ that minimises the number of blue edges in $X_1$ or $X_2$. Let $R_i$ and $B_i$ be the red and blue neighbourhoods of $u$ in $X_i$, respectively, for $i \in [2]$. Without loss of generality, $|B_1| \le |B_2|$.
				We consider two cases according to the size of $B_2$. 

				{\bf $|B_2| \ge n/10$}.\,\,
				Write $X_1' = X_1 \cup \{u\}$ and $X_2' = X_2$. 
				We claim that there is a triangle packing $\T_B$ that consists of blue cross triangles and that covers all blue edges in either $X_1'$ or $X_2'$.
				Indeed, by \Cref{prop:bip-stability-blue-neighs}~\ref{itm:stability-few-red-neighs}, there is a matching $M$ in $H_B[B_1, B_2]$ that covers $B_1$. Now, by \Cref{prop:bip-stability-blue-neighs}~\ref{itm:stability-perfect-blue-packing}, there is a triangle packing $\T_B'$ that consists of cross blue triangles, and covers all blue triangles in either $X_1$ or $X_2$. Add to $\T_B'$ the packing obtained by attaching $u$ to the edges of $M$, to obtain the required triangle packing $\T_B$    

				As usual, by \Cref{thm:n-four}, $G_R[X_i'] \setminus \T_B$ has a fractional triangle decompositions, for $i \in [2]$. It follows that 
				\begin{equation*} 
					\fpack(G) \ge (n+1)^2/4 - (n+1)/2 + 2|\T_B| = (n^2 - 1)/4 + 2k.
				\end{equation*}
				\Cref{lem:bip-stability} readily follows: if $k > \ell$, then $\fpack(G) \ge (n^2 - 1)/4 + 2\ell$; and if $k \le \ell$, then $G_B$ is $\ell$-close to bipartite.

				{\bf $|B_2| \le n/10$}.\,\,
				Let $\T_B$ be a triangle packing in $H$ that consists of blue cross triangles, avoids the edges of $M$, and covers all blue edges in either $X_1$ or $X_2$; such a packing exists due to \Cref{prop:bip-stability-blue-neighs}~\ref{itm:stability-perfect-blue-packing}.
				Let $M_i$ be a largest matching in $H_R[R_i]$, for $i \in [2]$. Then $M_i$ covers at least $|R_i| - \sqrt{n}$ vertices. 
				Note that $G_R[X_i] \setminus M_i$ has a triangle decomposition, by \Cref{thm:n-four}. 
				Consider the red triangle packing obtained by attaching $u$ to the edges of $M_1 \cup M_2$. It follows that $G$ has a red fractional triangle packing of size at least 
				\begin{align*} 
					n^2/4 - n/2 + 2k + 2(|M_1| + |M_2|)
					& \ge n^2/4 - n/2 + |R_1| + |R_2| - 2\sqrt{n} \\
					& \ge (n^2 - 1)/4 - n/2 + 4n/5 - 2\sqrt{n}, \\
					& = (n^2 - 1)/4 + 3n/10 - 2\sqrt{n} \\
					& > (n^2 - 1)/4 + 2(n/8 - n/200).
				\end{align*}
				It follows that $\fpack(G) \ge (n^2 - 1)/4 + 2\ell$ for any $\ell \le (1/8 + 1/200)n$, as required.
			\end{proof}

	\section{Pentagon blow-ups} \label{sec:pentagon}

		Given a red-blue colouring of a complete graph $G$, we say that $G$ is a \emph{pentagon blow-up} if there is a partition $\{A_1, \ldots, A_5\}$ of $V(G)$ into non-empty sets such that $G_R[A_i, A_{i+1}]$ and $G_B[A_i, A_{i+2}]$ are complete bipartite graphs, for $i \in [5]$ (where addition of indices is taken modulo $5$). Similarly, if $A_1, \ldots, A_5$ are non-empty pairwise disjoint sets in a red-blue colouring $G$ of a complete graph, we say that $(A_1, \ldots, A_5)$ is a \emph{pentagon blow-up} if the aforementioned property holds.
		Define
		\begin{align*}
			\B_1 =  
			\{& (3,3,3,4,4),\, (2,3,4,4,4),\, (3,3,3,3,5),\, (3,3,4,4,4),\, (2,4,4,4,4),\, (3,3,3,4,5),\, \\
			&(3,4,4,4,4),\, (3,3,4,4,5),\, (4,4,4,4,4),\, (3,4,4,4,5),\, (4,4,4,4,5),\, (3,4,4,5,5),\, \\
			&(4,4,4,5,5),\, (4,4,5,5,5),\, (4,5,5,5,5),\, (5,5,5,5,5)\} \\
			\B_2 = 
			\{& (3,3,3,4,4),\, (3,3,4,4,4),\, (3,4,4,4,4),\, (4,4,4,4,4),\, (4,4,4,4,5)\}.
		\end{align*}
		Let $\C_1$ be the family of pentagon blow-ups with blob sizes $x_1, \ldots, x_5$ (in some order), where $(x_1, \ldots, x_5) \in \B_1$, and let $\C_2$ be the family of red-blue colourings of $K_n$ that are one edge-flip away from a pentagon blow-up with blob sizes $x_1, \ldots, x_5$ (in some order), where $(x_1, \ldots, x_5) \in \B_2$. Set $\C = \C_1 \cup \C_2$.

		The main result in this section is the following theorem. 

		\begin{thm} \label{thm:pentagon}
			Let $G$ be a red-blue colouring of $K_n$, let $u \in V(G)$, and set $H = G \setminus \{u\}$. If $H \in \C$ then either $\fpack(G) > n(n-1)/4$ or $G \in \C$.
		\end{thm}

		Note that \Cref{lem:pentagon} follows from \Cref{thm:pentagon}, as the families in $\C$ have order at most $25$.

		\subsection{Fractional triangle packings between two or there blobs}

			Naturally, we would like to calculate the maximum size of a monochromatic fractional triangle packing of families in $\C$ (and related graphs). The following two propositions, regarding fractional triangle packings between two and three blobs, will allow us to do so. While the statement of the following proposition is a bit more general, we are really only interested in a finite number of small graphs (with blob sizes up to $6$), which we could have dealt with by computer search. Nevertheless, we include human-readable proofs of the following two propositions in \Cref{sec:two-three-blobs}.

			\begin{prop} \label{prop:AB}
				Let $G$ be a graph with vertex set $A \cup B$, where $A$ and $B$ are disjoint. In each of the following cases there is a fractional triangle packing $\omega$ in $G$, consisting of \emph{cross triangles} (namely, triangles with a vertex in both $A$ and $B$), such that $\omega(e) = 1/2$ for every edge $e$ in $A$ or in $B$.
				\begin{enumerate}[label = \rm(\alph*)]
					\item \label{itm:AB-complete}
						$2 \le |A| \le |B| \le |A| + 2$, and $G[A, B]$ is complete bipartite.
					\item \label{itm:AB-balanced-matching}
						$3 \le |A| \le |B| \le |A| + 1$, and $G[A, B]$ is complete bipartite minus a matching.
					\item \label{itm:AB-balanced-two-edges}
						$3 \le |A| \le |B| \le |A| + 1$, and $G[A, B]$ is complete bipartite minus two edges that intersect at $A$.
					\item \label{itm:AB-three-five}
						$|A| = 3$, $|B| = 5$, and $G[A, B]$ is complete bipartite minus a matching of size $2$.
				\end{enumerate}
			\end{prop}

			\begin{prop} \label{prop:ABC}
				Let $G$ be a graph with vertex set $A \cup B \cup C$, where $A, B, C$ are disjoint sets, with $|A| = 2$, $|B|, |C| \in \{3, 4\}$. Suppose that $G[B, A \cup C]$ is complete bipartite minus a matching with at most two edges between $B$ and $C$, and $G[A, C]$ is empty. Then there is a fractional triangle packing $\omega$ in $G$, consisting of cross triangles (namely, triangles intersecting (exactly) two of the sets $A, B, C$), such that $\omega(e) = 1/2$ if $e$ is in $A$ or in $C$, and $\omega(e) = 1$ if $e$ is in $B$.
			\end{prop}

			We use \Cref{prop:AB,prop:ABC} to calculate $\fpack(G)$ for every $G \in \C$; these values also appear in \Cref{table:C}.

			\begin{prop} \label{prop:pentagon-packing}
				\hfill
				\begin{enumerate}[label = \rm(\alph*)]
					\item \label{itm:pentagon-packing}
						Let $G \in \C_1$, so $G$ is a pentagon blow-up with blob sizes $x_1, \ldots, x_5$, where $(x_1, \ldots, x_5) \in \B_1$. Then $\fpack(G) = 3\sum_{i \in [5]} \binom{x_i}{2}$.
					\item \label{itm:almost-pentagon-packing}
						Let $G \in \C_2$, so $G$ is one edge-flip away from a pentagon blow-up with blob sizes $x_1, \ldots, x_5$, where $(x_1, \ldots, x_5) \in \B_2$. Then $\fpack(G) = 3\left(\sum_{i \in [5]} \binom{x_i}{2} + 1\right)$.
				\end{enumerate}
			\end{prop}

			\begin{proof}[Proof of \Cref{prop:pentagon-packing}]
				Let $\{A_1, \ldots, A_5\}$ be a partition of $V(G)$ such that $G$ is a pentagon blow-up if \ref{itm:pentagon-packing} holds, or one edge-flip away from a pentagon blow-up if \ref{itm:almost-pentagon-packing} holds, with blobs $A_1, \ldots, A_5$. So $A_1, \ldots, A_5$ have sizes $x_1, \ldots, x_5$ (possibly in a different order) where $(x_1, \ldots, x_5) \in \B_1$ if \ref{itm:pentagon-packing} holds, and $(x_1, \ldots, x_5) \in \B_2$ if \ref{itm:almost-pentagon-packing} holds. Also, $G_R[A_i, A_{i+1}]$ and $G_B[A_i, A_{i+2}]$ are complete bipartite if \ref{itm:pentagon-packing} holds; and if \ref{itm:almost-pentagon-packing} holds, without loss of generality there exist $x \in A_5$ and $y \in A_2$ such that $xy$ is red, and if $xy$ is recoloured blue then $G_R[A_i, A_{i+1}]$ and $G_B[A_i, A_{i+2}]$ are complete bipartite for $i \in [5]$. Let $z \in A_1$.

				Note that every monochromatic triangle in $G$ has an edge with both ends in $A_i$, for some $i \in [5]$, or contains $xy$ if \ref{itm:almost-pentagon-packing} holds. Thus $\fpack(G) \le 3\sum_{i \in [5]} \binom{x_i}{2}$ in case \ref{itm:pentagon-packing}; and $\fpack(G) \le 3\left(\sum_{i \in [5]} \binom{x_i}{2} + 1\right)$ in case \ref{itm:almost-pentagon-packing}. It remains to prove matching lower bounds. 
				
				If \ref{itm:pentagon-packing} holds, by \Cref{prop:AB}~\ref{itm:AB-complete}, there exist fractional triangle packings in $G_R[A_i, A_{i+1}]$ and $G_B[A_i, A_{i+2}]$ that consist of cross triangles and assign weight $1/2$ to each edge in $A_i$, $A_{i+1}$ or $A_{i+2}$. Putting these ten packings together, we obtain a monochromatic fractional triangle packing in $G$, that consists of cross triangles (i.e.\ triangles with vertices in two sets $A_i$) and assigns weight $1$ to each edge in $A_i$ for $i \in [5]$. It follows that $\fpack(G) \ge 3\sum_{i \in [5]}\binom{x_i}{2}$, as required for \ref{itm:pentagon-packing}.

				If \ref{itm:almost-pentagon-packing} holds, by \Cref{prop:AB}~\ref{itm:AB-balanced-matching}, there exist fractional triangle packings in $G_R[A_i, A_{i+1}]$ and $G_B[A_i, A_{i+2}]$ that consist of cross triangles that do not contain the edges $xy$, $xz$ or $yz$, and assign weight $1/2$ to each edge in $A_i$, $A_{i+1}$ or $A_{i+2}$. Putting together these packings, and adding the triangle $xyz$ with weight $1$, we find that $\fpack(G) \ge 3\left(\sum_{i \in [5]} \binom{x_i}{2} + 1\right)$, as required for \ref{itm:almost-pentagon-packing}.
			\end{proof}

		\subsection{Extensions of pentagon blow-ups and almost pentagon blow-ups}

			A \emph{balanced pentagon blow-up} is a pentagon blow-up whose blob sizes differ by at most one.
			\Cref{lem:pentagon-balanced,lem:pentagon-two-to-five,lem:pentagon-three-four-four-four-five,lem:pentagon-balanced-mistake} below will allow us to prove \Cref{thm:pentagon} for different families in $\C$. 
			The proofs of these lemmas are very similar (though the proof of \Cref{lem:pentagon-balanced-mistake} is a little more complicated). We prove \Cref{lem:pentagon-balanced} in \Cref{subsec:balanced-pentagon-blowup}, after deducing \Cref{thm:pentagon} from \Cref{lem:pentagon-balanced,lem:pentagon-two-to-five,lem:pentagon-three-four-four-four-five,lem:pentagon-balanced-mistake} and introducing some preliminaries. We delay the proofs of \Cref{lem:pentagon-two-to-five,lem:pentagon-three-four-four-four-five,lem:pentagon-balanced-mistake} to \Cref{appendix:extend-pentagon-blowups}.

			\begin{lem} \label{lem:pentagon-balanced}
				Let $n \ge 15$, let $G$ be a red-blue colouring of $K_n$, let $u \in V(G)$, and set $H = G \setminus \{u\}$.
				Suppose that $H$ is a balanced pentagon blow-up with blobs $A_1, \ldots, A_5$. Then the following holds, where $t = \floor{n/5}$. 
				\begin{itemize}
					\item
						$\fpack(G) \ge \fpack(H) + 3t$.
					\item \label{itm:t+one}
						$\fpack(G) \ge \fpack(H) + 3(t + 1)$, unless $(A_i \cup \{u\}, A_{i+1}, \ldots, A_{i+4})$ is a pentagon blow-up for some $i$ with $|A_i| = t$.
					\item \label{itm:t+two}
						$\fpack(G) \ge \fpack(H) + 3(t + 2)$, unless $(A_i \cup \{u\}, A_{i+1}, \ldots, A_{i+4})$ is one edge-flip away from a pentagon blow-up for some $i$ with $|A_i| = t$, or $(A_i \cup \{u\}, A_{i+1}, \ldots, A_{i+4})$ is a pentagon blow-up for some $i$.
				\end{itemize}
			\end{lem}

			\begin{lem} \label{lem:pentagon-two-to-five}
				Let $G$ be a red-blue colouring of $K_n$, let $u \in V(G)$, and set $H = G \setminus \{u\}$. Suppose that $H$ is a pentagon blow-up with blobs $A_1, \ldots, A_5$ of sizes in $\{3,4,5\}$, or in $\{2,3,4\}$ with at most one blob of size $2$. 
				Then $\fpack(G) \ge \fpack(H) + 12$, unless 
				\begin{itemize}
					\item
						$(A_i \cup \{u\}, A_{i+1}, \ldots, A_{i+4})$ is a pentagon blow-up for some $i$ with $|A_i| \le 3$, or 
					\item
						$(A_i \cup \{u\}, A_{i+1}, \ldots, A_{i+4})$ is one edge-flip from a pentagon blow-up for some $i$ with $|A_i| = 2$.
				\end{itemize}
			\end{lem}

			\begin{lem} \label{lem:pentagon-three-four-four-four-five}
				Let $G$ be a red-blue colouring of $K_n$, let $u \in V(G)$, and set $H = G \setminus \{u\}$. Suppose that $H$ is a pentagon blow-up with blobs $A_1, \ldots, A_5$ of sizes $3, 4, 4, 4, 5$ (not necessarily in this order). Then $\fpack(G) \ge \fpack(H) + 15$, unless
				\begin{itemize} 
					\item
						$(A_i \cup \{u\}, A_{i+1}, \ldots, A_{i+4})$ is a pentagon blow-up, for some $i$ with $|A_i| \le 4$, or
					\item
						$G$ is one edge-flip away from a balanced pentagon blow-up.
				\end{itemize}
			\end{lem}

			\begin{lem} \label{lem:pentagon-balanced-mistake}
				Let $n \ge 15$, let $G$ be a red-blue colouring of $K_n$, let $u \in V(G)$, and set $H = G \setminus \{u\}$.
				Suppose that $H$ is one edge-flip away from a balanced pentagon blow-up with blobs $A_1, \ldots, A_5$.
				Then the following holds, where $t = \floor{n/5}$.
				\begin{itemize}
					\item
						$\fpack(G) \ge \fpack(H) + 3t$,
					\item
						$\fpack(G) \ge \fpack(H) + 3(t + 1)$, unless $(A_i \cup \{u\}, A_{i+1}, \ldots, A_{i+4})$ is one edge-flip away from a pentagon blow-up for some $i$ with $|A_i| = t$. 
				\end{itemize}
			\end{lem}

		\subsection{Proof of \Cref{thm:pentagon}}

			\begin{proof}[Proof of \Cref{thm:pentagon}]
				\Cref{thm:pentagon} follows easily, but quite tediously, from \Cref{lem:pentagon-balanced,lem:pentagon-two-to-five,lem:pentagon-three-four-four-four-five,lem:pentagon-balanced-mistake}. To simplify the verification process, we included \Cref{table:C}. For example, one can read of it that any graph $H$ which is one edge-flip away from a pentagon blow-up with blob sizes $4,4,4,4,5$ cannot be extended to a red-blue colouring $G$ of $K_{22}$ with $\fpack(G) \le (22 \cdot 21)/4$, due to \Cref{lem:pentagon-balanced-mistake}.
			\end{proof}

			\begin{table}[h!] 
				\centering
				\begin{tabular}{|c|c|c|c|c|} 
					\hline
					$n$ & graph & $\fpack(\cdot)$ & $\frac{n(n+1)}{4}$ & $\floor{\frac{(n-1)^2}{4}}$ \\
					\hline
					\hline
					\multirow{4}{*}{$17$} & \cellcolor{\cola} \hyperref[lem:pentagon-balanced]{$(33344)$}\phantom{$*$} & $63$ & \multirow{4}{*}{$76.5$} & \multirow{4}{*}{$64$} \\
					& \cellcolor{\colb} \hyperref[lem:pentagon-balanced-mistake]{$(33344)*$} & $66$ & & \\
					& \cellcolor{\colc} \hyperref[lem:pentagon-two-to-five]{$(23444)$}\phantom{$*$} & $66$ & & \\
					& \cellcolor{\colc} \hyperref[lem:pentagon-two-to-five]{$(33335)$}\phantom{$*$} & $66$ & & \\
					\hline
					\multirow{4}{*}{$18$} & \cellcolor{\cola} \hyperref[lem:pentagon-balanced]{$(33444)$}\phantom{$*$} & $72$ & \multirow{4}{*}{$85.5$} & \multirow{4}{*}{$72$} \\
					& \cellcolor{\colb} \hyperref[lem:pentagon-balanced-mistake]{$(33444)*$} & $75$ & & \\
					& \cellcolor{\colc} \hyperref[lem:pentagon-two-to-five]{$(24444)$}\phantom{$*$} & $75$ & & \\
					& \cellcolor{\colc} \hyperref[lem:pentagon-two-to-five]{$(33345)$}\phantom{$*$} & $75$ & & \\
					\hline
					\multirow{3}{*}{$19$} & \cellcolor{\cola} \hyperref[lem:pentagon-balanced]{$(34444)$}\phantom{$*$} & $81$ & \multirow{3}{*}{$95$}& \multirow{3}{*}{$81$} \\
					& \cellcolor{\colb} \hyperref[lem:pentagon-balanced-mistake]{$(34444)*$} & $84$ & & \\
					& \cellcolor{\colc} \hyperref[lem:pentagon-two-to-five]{$(33445)$}\phantom{$*$} & $84$ & & \\
					\hline
					\multirow{3}{*}{$20$} & \cellcolor{\cola} \hyperref[lem:pentagon-balanced]{$(44444)$}\phantom{$*$} & $90$ & \multirow{3}{*}{$105$} & \multirow{3}{*}{$90$} \\
					& \cellcolor{\colb} \hyperref[lem:pentagon-balanced-mistake]{$(44444)*$} & $93$ & & \\
					& \cellcolor{\cold} \hyperref[lem:pentagon-three-four-four-four-five]{$(34445)$}\phantom{$*$} & $93$ & & \\
					\hline
					\multirow{3}{*}{$21$} & \cellcolor{\cola} \hyperref[lem:pentagon-balanced]{$(44445)$}\phantom{$*$} & $102$ & \multirow{3}{*}{$115.5$} & \multirow{3}{*}{$100$} \\
					& \cellcolor{\colb} \hyperref[lem:pentagon-balanced-mistake]{$(44445)*$} & $105$ & & \\
					& \cellcolor{\colc} \hyperref[lem:pentagon-two-to-five]{$(34455)$}\phantom{$*$} & $105$ & & \\
					\hline
					$22$ & \cellcolor{\cola} \hyperref[lem:pentagon-balanced]{$(44455)$}\phantom{$*$} & $114$ & $126.5$ & $110$ \\
					\hline
					$23$ & \cellcolor{\cola} \hyperref[lem:pentagon-balanced]{$(44555)$}\phantom{$*$} & $126$ & $138$ & $121$ \\
					\hline
					$24$ & \cellcolor{\cola} \hyperref[lem:pentagon-balanced]{$(45555)$}\phantom{$*$} & $138$ & $150$ & $132$ \\
					\hline
					$25$ &\cellcolor{\cola} \hyperref[lem:pentagon-balanced]{$(55555)$}\phantom{$*$} & $150$ & $162.5$ & $144$ \\
					\hline
				\end{tabular}
				\caption{A table depicting the families of graphs in $\C$ along with the size of their largest monochromatic fractional triangle packing, their number of vertices $n$, and the values $n(n+1)/4$ and $\floor{(n-1)^2/4}$. 
				\\ \\ 
				Here $(x_1 \ldots x_5)$ refers to the family of pentagon blow-ups with blobs sizes $x_1, \ldots, x_5$ (in some order), and $(x_1 \ldots x_5)*$ denotes the family of graphs that are one edge-flip away from a pentagon blow-up with blob sizes $x_1, \ldots, x_5$. 
				\\ \\ 
				The colours of a cell points to the lemma relevant to the family the cell represents: yellow points to \Cref{lem:pentagon-balanced}, red points to \Cref{lem:pentagon-two-to-five}, blue points to \Cref{lem:pentagon-three-four-four-four-five}, and green points to \Cref{lem:pentagon-balanced-mistake}.}
				\label{table:C}
			\end{table}

		\subsection{Bad configurations}

			In this subsection we introduce the notion of \emph{bad configurations} which will be handy in the proof of \Cref{lem:pentagon-balanced,lem:pentagon-two-to-five,lem:pentagon-three-four-four-four-five,lem:pentagon-balanced-mistake}.

			Let $G$ be a red-blue colouring of $K_n$, let $u \in V(G)$, and let $A_1, \ldots, A_5$ be pairwise disjoint non-empty sets in $V(G) \setminus \{u\}$. Suppose that $(A_1, \ldots, A_5)$ is a pentagon blow-up. A \emph{bad configuration} in $(A_1, \ldots, A_5)$ with respect to $u$ is either a set of three red neighbours of $u$,one from each set $A_i, A_{i+2}, A_{i+3}$, for some $i \in [5]$, or a set of three blue neighbours of $u$, one from each set $A_i, A_{i+1}, A_{i+2}$, for some $i \in [5]$.

			\begin{obs} \label{obs:quintuple}
				Let $G$ be a red-blue colouring of $K_n$, let $u \in V(G)$ and let $a_1, \ldots, a_5$ be distinct vertices in $V(G) \setminus \{u\}$ that form a pentagon blow-up. Then $\{u, a_1, \ldots, a_5\}$ spans a monochromatic triangle (which contains $u$). Moreover, if there is a bad configuration among $a_1, \ldots, a_5$ (with respect to $u$), then there are two edge-disjoint monochromatic triangles in $\{u, a_1, \ldots, a_5\}$ (which contain $u$).
			\end{obs}

			\begin{proof} 
				We assume that $a_i a_{i+1}$ is red for every $i \in [5]$, so $a_i a_{i+2}$ is blue for every $i \in [5]$.
				Without loss of generality, $ua_2$ is red. If $ua_1$ is red then $ua_1a_2$ is a red triangle, and similarly if $ua_3$ is red then $ua_2a_3$ is a red triangle. Otherwise, $ua_1$ and $ua_3$ are blue, so $ua_1a_3$ is a blue triangle. Regardless, there is a monochromatic triangle touching $u$.

				Now suppose that $(a_1, \ldots, a_5)$ contain a bad configuration. Without loss of generality, $a_2, a_4, a_5$ are red neighbours of $u$. Then $ua_4a_5$ is a red triangle, and, as above, $\{u, a_1, a_2, a_3\}$ span a monochromatic triangle that contains $u$ (and which is edge-disjoint of $ua_4a_5$). 
			\end{proof}

			\begin{prop} \label{prop:no-bad-configs}
				Let $G$ be a red-blue colouring of $K_n$, let $u \in V(G)$, and let $A_1, \ldots, A_5$ be pairwise disjoint sets in $V(G) \setminus \{u\}$ that form a pentagon blow-up in $G$. Suppose that $(A_1, \ldots, A_5)$ does not have a bad configuration with respect to $u$. Then $(A_i \cup \{u\}, A_{i+1}, \ldots, A_{i+4})$ is a pentagon blow-up for some $i$.
			\end{prop}

			\begin{proof} 
				Pick $a_i \in A_i$, for $i \in [5]$. Without loss of generality, at least three of $a_1, \ldots, a_5$ are red neighbours of $u$. In particular, $a_i$ and $a_{i+1}$ are red neighbours of $u$ for some $i \in [5]$; without loss of generality $i = 1$. If $a_4$ is a red neighbour of $u$, then $a_1, a_2, a_4$ form a bad configuration. Thus, without loss of generality, $a_3$ is a red neighbour of $u$. It follows from the lack of bad configurations that all vertices in $A_4 \cup A_5$ are blue neighbours in $u$, from which it follows that all vertices $A_1 \cup A_3$ are blue neighbours of $u$. We deduce that $(A_1, A_2 \cup \{u\}, A_3, A_4, A_5)$ is a pentagon blow-up. 
			\end{proof}

			\begin{prop} \label{prop:exactly-t-bad-configs}
				Let $G$ be a red-blue colouring of $K_n$, let $u \in V(G)$, and let $A_1, \ldots, A_5$ be pairwise disjoint sets in $V(G) \setminus \{u\}$ that form a pentagon blow-up in $G$. Suppose that there are $t$ pairwise disjoint bad configurations, but no $t+1$ pairwise disjoint bad configurations, where $t < \min_{i \in [5]} |A_i|$. Then there exists $i \in [5]$ such that $(A_i \cup \{u\}, A_{i+1}, \ldots, A_{i+4})$ is exactly $t$ edge-flips away from a pentagon blow-up.
			\end{prop}

			\begin{proof} 
				Let $S_1, \ldots, S_t$ be pairwise disjoint bad configurations, and denote $A_i' = A_i \setminus (S_1 \cup \ldots \cup S_t)$. By \Cref{prop:no-bad-configs}, there exists $i \in [5]$ such that $(A_i' \cup \{u\}, A_{i+1}', \ldots, A_{i+4}')$ is a pentagon blow-up; without loss of generality $i = 1$.
				We claim that for every $j \in [t]$, if $S_j$ consists of three red neighbours of $u$ then exactly one of them is in $A_3 \cup A_4$, and if $S_j$ consists of three blue neighbours of $u$ then exactly one of them is in $A_2 \cup A_5$. It would follow from this that $(A_1 \cup \{u\}, A_2, \ldots, A_5)$ is exactly $t$ edge-flips away from a pentagon blow-up. Indeed, each bad configuration gives rise to exactly one such edge-flip.

				Fix some $j \in [t]$. Without loss of generality, $S_j = \{a_{\ell}, a_{\ell+2}, a_{\ell+3}\}$, where $ua_i$ is red and $a_i \in A_i$ for $i \in \{\ell, \ell+2, \ell+3\}$ for some $\ell \in [5]$. Note that $S_j$ intersects $A_3 \cup A_4$. It remains to show that $S_j$ contains no more than one vertex from $A_3 \cup A_4$; or, equivalently, that $\ell \neq 1$. So suppose that $\ell = 1$. Then $|A_2'|, |A_5'| \ge 2$ (as $|A_2|, |A_5| \ge t+1$, and $S_1 \cup \ldots \cup S_t$ intersects each of $A_2$ and $A_5$ in at most $t-1$ vertices). Let $a_2, a_2' \in A_2$ and $a_5, a_5' \in A_5$ be distinct. Then $a_2, a_2', a_5, a_5'$ are red neighbours of $u$, and so $\{a_2, a_3, a_5\}$ and $\{a_2', a_4, a_5'\}$ are two disjoint bad configurations that are disjoint of $S_1, \ldots, S_{j-1}, S_{j+1}, \ldots, S_t$, implying that there are $t+1$ pairwise disjoint bad configurations, a contradiction.
			\end{proof}

		\subsection{Extending a balanced pentagon blow-up} \label{subsec:balanced-pentagon-blowup}

			\begin{proof}[Proof of \Cref{lem:pentagon-balanced}] 
				Note that if $(A_i \cup \{u\}, A_{i+1}, \ldots, A_{i+4})$ is a pentagon blow-up for some $i$ then, by \Cref{prop:AB}~\ref{itm:AB-complete}, 
				\begin{equation*}
					\fpack(G)
					= 3\left(\binom{|A_i| + 1}{2} + \sum_{j \in [4]} \binom{|A_{i+j}|}{2} \right)
					= \fpack(H) + 3|A_i|.
				\end{equation*}
				If $(A_i \cup \{u\}, A_{i+1}, \ldots, A_{i+4})$ is one edge-flip away from a pentagon blow-up for some $i$ with $|A_i| = t$, then by \Cref{prop:AB}~\ref{itm:AB-balanced-matching} and \Cref{prop:pentagon-packing},
				\begin{equation*}
					\fpack(G) 
					= 3\left(\binom{|A_i|+1}{2} + \sum_{j \in [4]} \binom{|A_{i+j}|}{2} + 1\right)
					= \fpack(H) + 3(|A_i| + 1).
				\end{equation*}
				Hence, if one of the conclusions given in the statement holds, then \Cref{lem:pentagon-balanced} holds. We may thus assume that they do not hold. 
				We consider three cases regarding the number of bad configurations in $(A_1, \ldots, A_5)$ with respect to $u$. In each of these cases we find a monochromatic triangle packing $\T$ that consists of $t+2$ triangles, each of which contains $u$ and two vertices from different sets $A_i$. It would follow from \Cref{prop:AB}~\ref{itm:AB-balanced-matching} and \Cref{prop:pentagon-packing} that
				\begin{equation*} 
					\fpack(G) 
					\ge \fpack(H \setminus \T) + 3|\T| 
					= \fpack(H) + 3(t+2),
				\end{equation*}
				as required.

				{\bf Two disjoint bad configurations.}\,
					Let $S_1, \ldots, S_t$ be pairwise disjoint transversals of $A_1, \ldots, A_5$, such that $S_1$ and $S_2$ contain bad configurations. 
					Take $\T$ to be a monochromatic triangle packings that consists of two triangles in $S_i \cup \{u\}$ for $i \in [2]$, and one triangle in $S_i \cup \{u\}$ for $i \in \{3, \ldots, t\}$; such triangles exist due to \Cref{obs:quintuple}, and they all contain $u$.

				{\bf One bad configuration, no two disjoint ones.}\,
					By \Cref{prop:exactly-t-bad-configs}, without loss of generality, $(A_1 \cup \{u\}, A_2, \ldots, A_5)$ is exactly one edge-flip away from a pentagon blow-up; by our assumptions, $|A_1| = t+1$.
					Assume, without loss of generality, that $x \in A_3$ is a red neighbour of $u$, and let $y \in A_2$. 
					Let $v_1, \ldots, v_{t+1}$ be an enumeration of $A_1$, and let $w_1, \ldots, w_{t+1}$ be distinct vertices such that $w_i \in (A_2 \cup A_5) \setminus \{y\}$ if $uv_i$ is red $w_i \in (A_3 \cup A_4) \setminus \{x\}$ if $uv_i$ is blue. Let $\T$ be the monochromatic triangle packing $\{uxy, uv_1w_1, \ldots, uv_{t+1}w_{t+1}\}$.

				{\bf No bad configurations.}\,
					By \Cref{prop:exactly-t-bad-configs}, without loss of generality, $(A_1 \cup \{u\}, A_2, \ldots, A_5)$ is a pentagon blow-up, a contradiction.
			\end{proof}

	\section{Conclusion} \label{sec:conc}

		We proved that in every $2$-coloured $K_n$ there is a collection of $n^2/12 + o(n^2)$ edge-disjoint triangles. It is natural to ask for the exact minimum.
		\begin{qn} 
			What is the minimum number of pairwise edge-disjoint monochromatic triangles that are guaranteed to exist in every $2$-colouring of $K_n$?
		\end{qn}
		It seems plausible that minimisers have one colour class which is an almost balanced complete bipartite graph minus a matching, similarly to \Cref{thm:frac-main-extremal}. However, divisibility considerations are likely to come into play. For example, if $n = 4m+2$ for some integer $m$, then it is better to have the blue edges span an almost balanced complete bipartite graph $K_{2m, 2m+2}$, than to have them span a balanced complete bipartite graph $K_{2m+1, 2m+1}$, as in the former example each vertex is incident with an odd number of red edges, at least one of which remains uncovered, leaving a total of at least $2m+1$ edges from the two red blobs uncovered in every monochromatic triangle packing. 

		Erd\H{o}s \cite{erdos} also considered a similar question\footnote{Like \Cref{conj:erdos}, Erd\H{o}s attributes the question to Faudree, Ordman and himself.}: how many edge-disjoint monochromatic triangles \emph{of the same colour} are guaranteed to exist in every $2$-colouring of $K_n$? Erd\H{o}s believed that the answer should exceed $(1 + \eps)n^2/24$. This is indeed the case: it readily follows from our stability result, \Cref{thm:main-stability}. Jacobson conjectured (see \cite{erdos-et-al}) that the answer should be $n^2/20 + o(n^2)$. 
		\begin{conj}[see Conjecture 2 in \cite{erdos-et-al}]
			In every $2$-colouring of $K_n$ there is a collection of $n^2/20 + O(n^2)$ pairwise edge-disjoint monochromatic triangles of the same colour.
		\end{conj}
		This conjecture is asymptotically tight, as seen by considering a balanced pentagon blow-up where the number of red and blue edges within blobs is roughly the same.

		Of course, it would also be interesting, but possibly very challenging, to find the minimum number of edge-disjoint monochromatic copies of $H$ among $2$-coloured $K_n$. It may also make sense to consider the same question for $r$ colours.
		\begin{qn} 
			Given a graph $H$ and an integer $r$, how many edge-disjoint monochromatic copies of $H$ is one guaranteed to find in an $r$-colouring of $K_n$?
		\end{qn}

		Finally, we mention a vaguely related question of Yuster \cite{yuster13}: how many edge-disjoint directed triangles are there guaranteed to be in a regular tournament on $n$ vertices (where $n$ is odd)? The best known bounds to date are due to Akaria and Yuster \cite{akaria-yuster}, who showed that the answer lies between $n^2/11.43 + o(n^2)$ and $n^2/9 + o(n^2)$. They conjectured that the upper bound, which is achieved by the tournament with vertices $[n]$ (for $n$ odd), and arcs $(i, i+j)$ for $i \in [n]$, $j \in [(n-1)/2]$, where addition is taken modulo $n$. We note that showing that this tournaments has $n^2/9 + o(n^2)$ edge-disjoint directed triangles requires some work.
		\begin{conj}[Conjecture 1.1 in \cite{akaria-yuster}]
			Every regular tournament on $n$ vertices, where $n$ is odd, contains a collection of $n^2/9 + o(n^2)$ pairwise edge-disjoint directed triangles.
		\end{conj}

	\bibliography{triangle-packing}

\providecommand{\bysame}{\leavevmode\hbox to3em{\hrulefill}\thinspace}
\providecommand{\MR}{\relax\ifhmode\unskip\space\fi MR }
\providecommand{\MRhref}[2]{%
  \href{http://www.ams.org/mathscinet-getitem?mr=#1}{#2}
}
\providecommand{\href}[2]{#2}
\begin{thebibliography}{10}

\bibitem{akaria-yuster}
I.~Akaria and R.~Yuster, \emph{Packing edge-disjoint triangles in regular and
  almost regular tournaments}, Discr. Math. \textbf{338} (2015), 217--228.

\bibitem{alon-shapira-sudakov}
N.~Alon, A.~Shapira, and B.~Sudakov, \emph{Additive approximation for
  edge-deletion problems}, Annals of Math. \textbf{170} (2009), 371--411.

\bibitem{barber-et-al}
B.~Barber, D.~K\"uhn, A.~Lo, and D.~Osthus, \emph{Edge-decompositions of graphs
  with high minimum degree}, Adv. Math. \textbf{288} (2016), 337--385.

\bibitem{delcourt-postle}
M.~Delcourt and L.~Postle, \emph{{Progress towards Nash-Williams' Conjecture on
  Triangle Decompositions}}, arXiv:1909.00514 (2019).

\bibitem{dross}
F.~Dross, \emph{Fractional triangle decompositions in graphs with large minimum
  degree}, SIAM J. Discr. Math. \textbf{30} (2015), 36--42.

\bibitem{dukes}
P.~Dukes, \emph{Rational decomposition of dense hypergraphs and some related
  eigenvalue estimates}, Linear Algebra and its applications \textbf{436}
  (2012), 3726--3746, (see arXiv:1108.1576 for an erratum).

\bibitem{erdos}
P.~Erd\H{o}s, \emph{Some recent problems and results in graph theory}, Discr.
  Math. \textbf{164} (1997), 81--85.

\bibitem{erdos-et-al}
P.~Erd\H{o}s, R.~J. Faudree, R.~J. Gould, M.~S. Jacobson, and J.~Lehel,
  \emph{Edge disjoint monochromatic triangles in $2$-colored graphs}, Discr.
  Math. \textbf{231} (2001), 135--141.

\bibitem{garaschuk}
K.~Garaschuk, \emph{Linear methods for rational triangle decompositions}, Ph.D.
  thesis, University of Victoria, 2014.

\bibitem{goodman}
A.~W. Goodman, \emph{On sets of acquaintances and strangers at any party},
  Amer. Math. Monthly \textbf{66} (1959), 778--783.

\bibitem{us2}
V.~Gruslys and S.~Letzter, \emph{Fractional triangle packings in almost
  complete graphs}, arXiv:2008.05313.

\bibitem{gustavsson}
T.~Gustavsson, \emph{Decompositions of large graphs and digraphs with high
  minimum degree}, Ph.D. thesis, University of Stockholm, 1991.

\bibitem{haxell-rodl}
P.~Haxell and V.~R\"odl, \emph{Integer and fractional packings in dense
  graphs}, Combinatorica \textbf{21} (2001), 13--38.

\bibitem{keevash-sudakov}
P.~Keevash and B.~Sudakov, \emph{Packing triangles in a graph and its
  complement}, J. Graph Theory \textbf{47} (2004), 203--216.

\bibitem{mckay-piperno}
B.~D. McKay and A.~Piperno, \emph{{Practical Graph Isomorphism, II}}, J.
  Symbolic Computation \textbf{60} (2013), 94--112.

\bibitem{nash-williams}
C.~S.~J. Nash-Williams, \emph{An unsolved problem concerning decomposition of
  graphs into triangles}, Combinatorial Theory and its Applications
  \textbf{III} (1070), 1179--1183.

\bibitem{tyomkyn}
M.~Tyomkyn, \emph{Many disjoint triangles in co-triangle-free graphs},
  arXiv:2001.00763 (2020).

\bibitem{yuster05}
R.~Yuster, \emph{{Asymptotically optimal $K_k$-packings of dense graphs via
  fractional $K_k$-decompositions}}, J. Combin. Theory Ser. B \textbf{95}
  (2005), 1--11.

\bibitem{yuster07}
\bysame, \emph{Packing cliques in graphs with independence number $2$}, Combin.
  Probab. Comput. \textbf{16} (2007), 805--817.

\bibitem{yuster13}
\bysame, \emph{Packing triangles in regular tournaments}, J. Graph Theory
  \textbf{74} (2013), 58--66.

\end{thebibliography}
	\bibliographystyle{amsplain}

	\appendix

	\section{Fractional triangle packings with two or three blobs} \label{sec:two-three-blobs}

		\begin{proof}[Proof of \Cref{prop:AB}]
			Throughout the proof, we may assume that $G[A]$ and $G[B]$ are complete. Write $\alpha = |A|$ and $\beta = |B|$.

			Let $G$ satisfy one of the properties \ref{itm:AB-complete} and \ref{itm:AB-balanced-matching}. We define $\omega$ as follows: for each cross triangle $T$, we define $\omega(T) = 1/2d$, where $d$ is the number of common neighbours of the two vertices of $T$ that belong to the same side. Clearly, $\omega(e) = 1/2$ for every edge $e$ in $A$ or $B$. It remains to check that $\omega(e) \le 1$ for every cross edge $e$. It will be useful to note that the total weight with respect to $\omega$ on cross edges is twice the weight on edges in $A$ or $B$, namely $\binom{\alpha}{2} + \binom{\beta}{2}$.

			For \ref{itm:AB-complete}, note that the cross edges all receive the same weight, which is 
			\begin{equation*}
				\frac{\binom{\alpha}{2} + \binom{\beta}{2}}{\alpha \beta}
				= \frac{\alpha(\alpha - 1) + \beta(\beta - 1)}{2\alpha \beta}.
			\end{equation*}
			If $\beta \le \alpha + 1$, then clearly the denominator is at most $2\alpha \beta$, implying that the expression is at most $1$, as required. If $\beta = \alpha + 2$, then the expression is
			\begin{equation*} 
				\frac{2\alpha^2 + 2\alpha + 2}{2(\alpha^2 + 2\alpha)}
				\le 1,
			\end{equation*}
			as required.

			For \ref{itm:AB-balanced-matching}, we assume that the non-edges in $G[A, B]$ form a matching that saturates $A$. If $|B| = |A|$, again all cross edges have the same weight, which is $\frac{\alpha(\alpha - 1)}{\alpha^2 - \alpha} = 1$.
			If $|B| = |A| + 1$, then the total weight on cross edges is $\binom{\alpha}{2} + \binom{\alpha+1}{2} = \alpha^2$, which is also the number of cross edges. Note that all cross edges that intersect two non-edges have the same weight, and all remaining cross edges (which are incident to the unique vertex in $B$ not incident with a non-edge) also have the same weight. As the average weight is $1$, it suffices to show that edges of the latter type have weight $1$. Indeed, as each such edge is in $2(\alpha - 1)$ cross triangles, each of which has weight $\frac{1}{2(\alpha - 1)}$, \ref{itm:AB-balanced-matching} follows.

			For \ref{itm:AB-balanced-two-edges}, denote the non-edges in $G[A, B]$ by $xy$ and $xz$ (so $x \in A$), and write $\alpha = |A|$ and $\beta = |B|$.
			Let $\omega_1$ be a fractional triangle packing in $G$, define as follows for every $a \in A \setminus \{x\}$, $b \in B \setminus \{y,z\}$ and $v \in \{y, z\}$.
			\begin{equation*} 
				\omega_1(xab) = \frac{1}{2(\beta - 2)}, \qquad \qquad \omega_1(avb) = \frac{1}{2(\alpha - 1)}.
			\end{equation*}
			We note that $\omega_1(e) = 1/2$ for every edge $e$ which is in $A$ and touches $x$, or which is in $B$ and touches $\{y, z\}$. Moreover, $\omega_1(f) \le 1$ for every edge $f$ between $A$ and $B$. Indeed, if $f = uw$ for $u \in A$ and $w \in B$, then
			\begin{equation*} 
				\omega_1(f) = \left\{
					\begin{array}{ll}
						\frac{\alpha - 1}{2(\beta - 2)} & u = x \\
						\frac{\beta - 1}{2(\alpha - 1)} & w \in \{y, z\} \\
						\frac{1}{2(\beta - 2)} + \frac{2}{2(\alpha - 1)} & \text{otherwise},
					\end{array}
				\right.
			\end{equation*}
			and one can easily check that $\omega_1(f) \le 1$ for every cross edge $f$, using $3 \le \alpha \le \beta \le \alpha + 1$.
			Next, let $\omega_2$ be the fractional triangle packing, defined by giving each triangle $xbb'$, where $b, b' \in B \setminus \{y, z\}$, the same weight, so that $\omega_1(f) + \omega_2(f) = 1$ for $f = xb$ with $b \in B \setminus \{y, z\}$ (in case $\beta = 3$, there are no such triangles and $\omega_2$ assigns weight $0$ to all triangles). Similarly, let $\omega_3$ be the fractional triangle packing, obtained by giving each triangle $vaa'$, where $v \in \{y, z\}$ and $a, a' \in A \setminus \{x\}$ the same weight, so that $\omega_1(f) + \omega_3(f) = 1$ for $f = va$ for $v \in \{y, z\}$ and $a \in A \setminus \{x\}$.
			As $\omega_1(f) \ge 1/2$ for $f = xb$ with $b \in B \setminus \{y, z\}$ or $f = va$ with $v \in \{y, z\}$ and $a \in A \setminus \{x\}$, we find that $\omega_2(e) \le 1/2$ for $e = bb'$ with $b, b' \in B \setminus \{y, z\}$, and $\omega_3(e) \le 1/2$ for $e = aa'$ for $a, a' \in A \setminus \{x\}$. Finally, we claim that there is a fractional triangle packing $\omega_4$, which assigns non-zero weight only to cross triangles that avoid $x, y, z$, so that $\omega = \omega_1 + \ldots + \omega_4$ is a fractional triangle packing with the required properties. 
			Indeed, it suffices to show, by symmetry, that the available weight on cross edges between $A \setminus \{x\}$ and $B \setminus \{y, z\}$ is at least twice the available weight in $A \setminus \{x\}$ and in $B \setminus \{y, z\}$. The latter statement follows, as $\omega_1, \omega_2, \omega_3$ saturate all edges that touch $\{x, y, z\}$ (with the exception of the single edge from $x$ to $B \setminus \{y, z\}$ if $\beta = 3$), and the total available weight on cross edges, minus $1$, is at least twice the total available weight in $A$ or in $B$ (we made the relevant calculation for \ref{itm:AB-balanced-matching}). We have thus established \ref{itm:AB-balanced-two-edges}.

			It remains to prove \ref{itm:AB-three-five}. For a cross triangle $T$, we define $\omega(T)$ as follows, where $A'$ and $B'$ are the vertices in $A$ and in $B$, respectively, that are incident with a non-edge. 
			\begin{equation*} 
				\omega(T) = \left\{
					\begin{array}{ll}
						1/2 & \text{$T$ has two vertices in $B'$ (and one in $A \setminus A'$)} \\
						1/3 & \text{$T$ has a vertex in $A'$, a vertex in $B'$ and another in $B \setminus B'$} \\
						0 & \text{$T$ has a vertex in $A'$, a vertex in $A \setminus A'$, and another in $B'$,} \\
						1/6 & \text{otherwise}.
					\end{array}
				\right.
			\end{equation*}
			One can check that every edge in $A$ or in $B$ receives weight exactly $1/2$, and every cross edge receives weight $1$. 
		\end{proof}

		\begin{proof}[Proof of \Cref{prop:ABC}]
			We assume, without loss of generality, that $G[A]$, $G[B]$ and $G[C]$ are complete. 
			We consider two cases: $|B| = 4$, and $|B| = 3$.

			In the former case, we may assume that both $G[A, B]$ and $G[B, C]$ have two missing edges. Let $B' = \{b_1, b_2\}$ be the set of vertices in $B$ that are incident with missing edges in $G[A, B]$. For a cross triangle $T$ in $G[A \cup B]$, define $\omega(T)$ as follows 
			\begin{equation*} 
				\omega(T) = \left\{
					\begin{array}{ll} 
						1/2 & \text{$T$ has one vertex in $A$, one in $B'$, and one in $B \setminus B'$} \\
						1/4 & \text{$T$ has one vertex in $A$, and two in $B \setminus B'$} \\
						1/4 & \text{$T$ has two vertices in $A$, and one in $B \setminus B'$}.
					\end{array}
				\right.
			\end{equation*}
			We note that $b_1b_2$ receives weight $0$, and all other edges in $A$ or in $B$ receive weight $1/2$. One can also check that the weight of every cross edge in $G[A, B]$ is $1$. Next, we consider $G[B, C]$. Let $c \in C$ be any vertex which is not incident with a missing edge (there are either one or two such vertices). By \Cref{prop:AB}~\ref{itm:AB-balanced-matching}, there exist fractional triangle packings $\omega_i$, for $i \in [2]$, that consist of cross triangles in $G[B, C] \setminus b_ic$, such that $\omega(e) = 1/2$ for every edge $e$ in $B$ or in $C$. Define, for a cross triangle $T$ in $G[B, C]$, 
			\begin{equation*} 
				\omega(T) = \left\{
					\begin{array}{ll}
						1/2 & T = b_1b_2c \\
						\frac{\omega_1(T)\, +\, \omega_2(T)}{2} & \text{otherwise}.
					\end{array}
				\right.
			\end{equation*}
			It is easy to check that $\omega$ satisfies the requirements.

			In the second case, we may assume that either there are two missing edges in $G[A, B]$ and one missing edge in $G[B, C]$; or there is one missing edge in $G[A, B]$ and two in $G[B, C]$. 
			If the former holds, we define $\omega(T) = 1/2$ for all cross triangles in $G[A \cup B]$ (there are three such triangles); this assigns weight $0$ to the edge in $B$ whose two ends are incident with missing edges in $G[A, B]$, and weight $1/2$ to all other edges in $A$ or in $B$. We can then define a fractional triangle packing in $G[B, C]$ as above, using \Cref{prop:AB}~\ref{itm:AB-balanced-matching}, to obtain the required packing.
			If the latter holds, we define $\omega(T) = 1/2$ for the two cross triangles in $G[A \cup B]$ that contain the vertex in $B$ which is incident to a non-edge in $G[A, B]$; and to the other cross triangles $T$ in $G[A \cup B]$ we assign weight $1/4$. It is easy to check that every edge in $A$ or in $B$ receives a weight of $1/2$, and that cross edges in $G[A, B]$ receive a weight of at most $1$. The required packing can thus be obtained by taking a packing as in \Cref{prop:AB}~\ref{itm:AB-balanced-matching} in $G[B \cup C]$.
		\end{proof}

	\section{Extending pentagon blow-ups} \label{appendix:extend-pentagon-blowups}

		\subsection{Proof of \Cref{lem:pentagon-two-to-five}}

			\begin{proof}
				Suppose that neither of the two conclusions holds.
				We consider three cases regarding the number of bad configurations in $(A_1, \ldots, A_5)$.
				In each of these cases, we find a monochromatic triangle packing $\T$ that consists of four triangles, each of which contains $u$ and two vertices from two different sets $A_i$, no three of which have edges between the same two sets $A_i$ and $A_j$. For such $\T$, it follows from \Cref{prop:AB}~\ref{itm:AB-balanced-matching} and \ref{itm:AB-three-five}, \Cref{prop:ABC} and \Cref{prop:pentagon-packing} that
				\begin{equation*} 
					\fpack(G) \ge \fpack(H \setminus \T) + 3|\T| = \fpack(H) + 12,
				\end{equation*}
				as required.

				{\bf Two disjoint bad configurations.}\,
				Let $S_1$ and $S_2$ be disjoint transversals of $A_1, \ldots, A_5$ that contain bad configurations. Take $\T$ to be a monochromatic triangle packing that consists of two triangles in each of $S_1 \cup \{u\}$ and $S_2 \cup \{u\}$; these exist by \Cref{obs:quintuple}, all the triangles chosen contain $u$, and no three of these triangles contain edges between the same two sets $A_i$ and $A_j$.

				{\bf One bad configuration, no two disjoint ones.}\,
				By \Cref{prop:exactly-t-bad-configs}, without loss of generality, $(A_1 \cup \{u\}, A_2, \ldots, A_5)$ is exactly one edge-flip away from a pentagon blow-up, and $ux$ is red for some $x \in A_3$. Let $y \in A_2$. 
				By our assumptions, $|A_1| \ge 3$. Let $v_1, v_2, v_3$ be three distinct elements from $A_1$, and let $w_1, w_2, w_3$ be distinct vertices such that $w_i \in (A_2 \cup A_5) \setminus \{y\}$ if $uv_i$ is red, $w_i \in (A_3 \cup A_4) \setminus \{x\}$ if $uv_i$ is blue, and each set $A_i$ contains at most two vertices $w_i$. The packing $\T = \{uxy, uv_1w_1, uv_2w_2, uv_3w_3\}$ satisfies the requirements. 

				{\bf No bad configurations.}\,
				By \Cref{prop:no-bad-configs}, without loss of generality, $(A_1 \cup \{u\}, A_2, \ldots, A_5)$ is a pentagon blow-up. By our assumptions, $|A_1| \ge 4$. Let $v_1, \ldots, v_4 \in A_1$ be distinct, and let $w_1, \ldots, w_4$ be distinct vertices such that $w_i \in A_2 \cup A_5$ if $uv_i$ is red and $w_i \in A_3 \cup A_4$ if $uv_i$ is blue, and every set $A_i$ contains at most two vertices $w_i$. Take $\T = \{uv_1w_1, \ldots, uv_4w_4\}$. 
			\end{proof}

		\subsection{Proof of \Cref{lem:pentagon-three-four-four-four-five}}

			\begin{proof}
				As before, we assume that neither of the two conclusions hold.
				In each of the following four cases, we find a monochromatic packing $\T$ that consists of five triangles that contain $u$ and two vertices from two different sets $A_i$, at most three of which contain an edge between the blobs of sizes $3$ and $5$. For such $\T$, it follows from \Cref{prop:AB}~\ref{itm:AB-balanced-matching} and \ref{itm:AB-three-five} and \Cref{prop:pentagon-packing} that
				\begin{equation*} 
					\fpack(G) \ge \fpack(H \setminus \T) + 3|\T| = \fpack(H) + 15,
				\end{equation*}
				as required.

				{\bf Three pairwise disjoint bad configurations.}\,
				Let $S_1, S_2, S_3$ be pairwise disjoint transversals in $A_1, \ldots, A_5$ that contain bad configurations. Take $\T'$ to be a monochromatic triangle packing that consists of two triangles in $S_i \cup \{u\}$ for each $i \in [3]$; these exist by \Cref{obs:quintuple} and all triangles in $\T'$ contain $u$. Note that at most three triangles in $\T'$ contain an edge between the blobs of sizes $3$ and $5$. Remove one such triangle (if one exists), to obtain the required triangle packing.

				{\bf Two pairwise disjoint bad configurations, no three disjoint ones.}\,
				Let $S_1$ and $S_2$ be disjoint transversals of $A_1, \ldots, A_5$ that contain bad configurations, and let $\T'$ be a monochromatic triangle packing that consists of two triangles in $S_i \cup \{u\}$ for each $i \in [2]$. Consider $A_i' = A_1 \setminus (S_1 \cup S_2)$. Then, by \Cref{prop:no-bad-configs}, without loss of generality, $(A_1' \cup \{u\}, A_2', \ldots, A_5')$ is a pentagon blow-up. Pick $v \in A_1'$. If $uv$ is red let $w \in A_2' \cup A_5'$, and, otherwise, let $w \in A_3' \cup A_4'$; if $A_1$ has size $3$ or $5$ we choose $w$ so that it does not belong to the other blob of size $3$ or $5$.
				$\T = \T' \cup \{uvw\}$ satisfies the requirements.

				{\bf One bad configuration, no two disjoint ones.}\, 
				By \Cref{prop:exactly-t-bad-configs}, without loss of generality $(A_1 \cup \{u\}, A_2, \ldots, A_5)$ is one edge-flip away from a pentagon blow-up, and $x \in A_3$ is a red neighbour of $u$. Let $y \in A_2$. By our assumptions, $|A_1| \ge 4$. Let $v_1, \ldots, v_4 \in A_1$ be distinct. Let $w_1, \ldots, w_4$ be distinct vertices such that $w_i \in (A_2 \cup A_5) \setminus \{y\}$ if $uv_i$ is red, $w_i \in (A_3 \cup A_4) \setminus \{x\}$ if $uv_i$ is blue, and if $|A_1| = 5$ then at most two vertices $w_i$ belong to the blob of size $3$. The triangle packing $\T = \{uxy, uv_1w_1, \ldots, uv_4w_4\}$ satisfies the requirements.

				{\bf No bad configurations.}\,
				By \Cref{prop:no-bad-configs}, without loss of generality, $(A_1 \cup \{u\}, A_2, \ldots, A_5)$ is a pentagon blow-up. By our assumptions, $|A_1| = 5$. Let $v_1, \ldots, v_5$ be an enumeration of $A_1$. Let $w_1, \ldots, w_5$ be disjoint vertices such that $w_i \in A_2 \cup A_5$ if $uv_i$ is red, $w_i \in A_3 \cup A_4$ if $uv_i$ is blue, and at most two vertices $w_i$ belong to the blob of size $3$. Take $\T = \{uv_1w_1, \ldots, uv_5w_5\}$. 
			\end{proof}

		\subsection{Proof of \Cref{lem:pentagon-balanced-mistake}}

			\begin{proof}
				Without loss of generality, there are vertices $x \in A_5$ and $y \in A_2$ such that $xy$ is red. 
				Suppose that the conclusion in the second item holds. Then, by \Cref{prop:pentagon-packing}, $\fpack(G) = \fpack(H) + 3t$, as required.

				We now assume that the conclusion does not hold. In each of the following three cases, we find a monochromatic triangle packing $\T$ that consists of $t+2$ triangles that intersect each blob $A_i$ in at most one vertex, one of which contains the edge $xy$, and such that $E(\T) \cap (A_i \times A_j)$ is either a matching or a set of two intersecting edges whose common vertex lies in a blob of size $t$. Then, by \Cref{prop:AB}~\ref{itm:AB-balanced-matching} and \ref{itm:AB-balanced-two-edges} and \Cref{prop:pentagon-packing},
				\begin{equation*} 
					\fpack(G) \ge \fpack(H \setminus \T) + 3|\T| \ge \fpack(H) + 3(t+1),
				\end{equation*}
				as required. 

				\subsubsection*{Two disjoint bad configurations in $(A_1', \ldots, A_5')$.}

					Let $S_1, \ldots, S_{t-1}$ be pairwise-disjoint transversals in $(A_1, A_2 \setminus \{y\}, A_3, A_4, A_5 \setminus \{x\})$, such that $S_1$ and $S_2$ contain bad configurations, and let $z \in A_1 \setminus (S_1 \cup \ldots \cup S_{t-1})$. Let $\T$ be a monochromatic triangle packing that consists of $xyz$ (a red triangle), two monochromatic triangles in $\{u\} \cup S_i$ for $i \in [2]$, and one monochromatic triangle in $\{u\} \cup S_i$ for $i \in \{3, \ldots, t-1\}$; such triangles exist due to \Cref{obs:quintuple}.

				\subsubsection*{One bad configuration, no two disjoint ones.}

					Let $S$ be a transversal in $(A_1, A_2 \setminus \{y\}, A_3, A_4, A_5 \setminus \{x\})$ that contains a bad configuration, let $A_i' = A_i \setminus (S \cup \{x, y\})$, and let $z \in A_1'$. Let $\T'$ be a triangle packing consisting of two edge-disjoint monochromatic triangles in $u \cup S$ that contain $u$; such $\T'$ exists by \Cref{obs:quintuple}. By \Cref{prop:no-bad-configs}, $(A_i' \cup \{u\}, A_{i+1}', \ldots, A_{i+4}')$ is a pentagon blow-up, for some $i \in [5]$. Without loss of generality, $i \in \{1, 2, 3\}$. 

					\begin{itemize} 
						\item
							$i = 3$.\,
							Let $v_1, \ldots, v_{t-1} \in A_3'$ be distinct, and let $w_1, \ldots, w_{t-1}$ be distinct vertices such that $w_j \in A_2' \cup A_4'$ if $uv_j$ is red, and $w_j \in A_1' \cup A_5'$ if $uv_j$ is blue.
							Set $\T = \{xyz, uv_1w_1, \ldots, uv_{t-1}w_{t-1}\} \cup \T'$.
						\item
							$i = 2$.\,
							Let $v_1, \ldots, v_{t-1} \in A_2' \cup \{y\}$ be distinct. Let $w_1, \ldots, w_{t-1}$ be distinct vertices such that $w_j \in A_3'$ if $uv_j$ is red and $w_j \in A_4'$ if $uv_j$ is blue. Take $\T = \{xyz, uv_1w_1, \ldots, uv_{t-1}w_{t-1}\} \cup \T'$. 
						\item
							$i = 1$.\,
							Without loss of generality, $\T'$ does not have a triangle with vertices in $A_1$ and $A_2$.
							Let $v_1, \ldots, v_{t-1} \in A_1'$ be distinct vertices such that $v_1, \ldots, v_{t-1} \neq z$ if $|A_1| = t+1$, and $v_{t-1} = z$ if $|A_1| = t$. Let $w_1, \ldots, w_{t-1}$ be distinct vertices such that 
							\begin{equation*} 
								\left\{
									\begin{array}{ll}
										w_j \in A_3' \cup A_4' & \text{$uv_j$ is blue and $j \in [t-1]$} \\
										w_j \in A_5' & \text{$uv_j$ is red and $j \in [t-2]$} \\
										w_j \in A_2' & \text{$uv_j$ is red and $j = t-1$}. 
									\end{array}
								\right.
							\end{equation*}
							Set $\T = \{xyz, uv_1w_1, \ldots, uv_{t-1}w_{t-1}\} \cup \T'$.
					\end{itemize}

				\subsubsection*{No bad configurations.}

					By \Cref{prop:no-bad-configs}, $(A_i' \cup \{u\}, A_{i+1}', \ldots, A_{i+4}')$ is a pentagon blow-up for some $i \in [5]$, where $A_i' = A_i \setminus \{x, y\}$. Again, without loss of generality, $i \in \{1,2,3\}$. 

					\begin{itemize}
						\item
							$i = 3$.\,
							Let $a \in A_5 \setminus \{x\}$, $b \in A_4$ and $z \in A_1$.
							Write $\ell = |A_3|$, and let $v_1, \ldots, v_{\ell}$ be an enumeration of the elements in $A_3$.   
							Let $w_1, \ldots, w_{\ell}$ be distinct vertices such that $w_j \in (A_2 \cup A_4) \setminus \{y, b\}$ if $uv_j$ is red, and  $w_j \in (A_1 \cup A_5) \setminus \{x, a\}$ if $uv_j$ is blue.
							Let $\T'$ be the monochromatic triangle packing $\{xyz, uv_1w_1, \ldots, uv_{\ell}w_{\ell}\}$. 
							By our assumptions, either $\ell = |A_3| = t+1$; $ux$ is red; or $uy$ is blue. If $|A_3| = t+1$ we take $\T = \T'$; if $ux$ is red we take $\T = \T' \cup \{uxb\}$; and if $uy$ is blue and $|A_2| = t$ we take $\T = \T' \cup \{uya\}$. It remains to consider the case where $|A_3| = t$, $ux$ is blue, $uy$ is blue, and $|A_2| = t+1$. In this case $G \setminus \{y\}$ is a balanced pentagon blow-up (with blobs $A_1, A_2 \setminus \{y\}, A_3 \cup \{u\}, A_4, A_5)$. The proof now follows from \Cref{lem:pentagon-balanced}.

						\item
							$i = 2$.\,
							Let $a, z \in A_1$ be distinct.
							Let $v_1, \ldots, v_t \in A_2$ be distinct.  Let $w_1, \ldots, w_t$ be distinct vertices such that $w_j \in A_3$ if $uv_j$ is red, and $w_j \in A_4$ if $uv_j$ is blue. Define $\T' = \{xyz, uv_1w_1, \ldots, uv_tw_t\}$. 
							By our assumptions, either $|A_2| = t+1$ or $ux$ is red.
							If the latter holds, take $\T = \T' \cup \{uxa\}$. Otherwise, $G \setminus \{y\}$ is a balanced pentagon blow-up, and the proof follows from \Cref{lem:pentagon-balanced}.
						
						\item
							$i = 1$.\,
							By our assumptions, either at least one of $ux$ and $uy$ is blue, or $|A_1| = t+1$.

							Suppose that the former holds; without loss of generality, $ux$ is blue.
							Let $a \in A_3$, $z \in A_1$. Let $v_1, \ldots, v_t \in A_1$ be distinct vertices such that $v_1, \ldots, v_t \neq z$ if $|A_1| = t+1$, and $v_t = z$ otherwise, and let $w_1, \ldots, w_t$ be distinct vertices such that 
							\begin{equation*} 
								\left\{
									\begin{array}{ll}
										w_j \in (A_3 \cup A_4) \setminus \{a\} & \text{$uv_j$ is blue and $j \in [t]$} \\
										w_j \in A_5 \setminus \{x\} & \text{$uv_j$ is red and $j \in [t-1]$} \\
										w_j \in A_2 \setminus \{y\} & \text{$uv_j$ is red and $j = t$}.
									\end{array}
								\right.
							\end{equation*}
							Take $\T = \{xyz, uxa, uv_1w_1, \ldots, uv_tw_t\}$.

							Now suppose that $ux$ and $uy$ are red and that $|A_1| = t+1$. Let $v_1, \ldots, v_{t+1}$ be an enumeration of the vertices in $A_1$, and let $w_1, \ldots, w_{t+1}$ be distinct vertices such that $w_j \in (A_2 \cup A_5) \setminus \{x,y\}$ if $uv_j$ is red, and $w_j \in A_3 \cup A_4$ if $uv_j$ is blue. Take $\T = \{uxy, uv_1w_1, \ldots, uv_{t+1}w_{t+1}\}$. 
							\qedhere
					\end{itemize}
			\end{proof}

\end{document}